\documentclass{amsart}
\usepackage{amssymb} 
\usepackage{color}

 \newtheorem{theorem}{Theorem}[section]

 \newtheorem{proposition}[theorem]{Proposition}
 \newtheorem{lemma}[theorem]{Lemma}
 \newtheorem{corollary}[theorem]{Corollary}
 \newtheorem{remark}[theorem]{Remark}
 \newtheorem{definition}[theorem]{Definition}
 \newtheorem{notation}[theorem]{Notation}
 \newtheorem{overview}[theorem]{Overview}
 \newtheorem{example}[theorem]{Example}

 \numberwithin{equation}{section}
 \def\subrel#1#2{\mathrel{\mathop{#2}\limits_{#1}}}
 \def\qqq{\qed\vspace{2mm}} 
 \def\power#1{[\![#1]\!]}

\begin{document}

\title{Hecke operators in Morava $E$-theories of different heights}
\author{Takeshi Torii}
\address{Department of Mathematics, 
Okayama University,
Okayama 700--8530, Japan}
\email{torii@math.okayama-u.ac.jp}
\thanks{This work was partially supported by 
JSPS KAKENHI Grant Number 22540087 and 25400092.}


\date{October 12, 2022\ ({\tt version~2.1})}

\begin{abstract}
There is a natural action
of a kind of Hecke algebra $\mathcal{H}_n$
on the $n$th Morava $E$-theory of spaces.
We construct Hecke operators 
in an amalgamated cohomology theory 
of the $n$th and the $(n+1)$st Morava $E$-theories.
These operations are natural extensions of the Hecke operators
in the $(n+1)$st Morava $E$-theory,
and they induce an action of the Hecke algebra 
$\mathcal{H}_{n+1}$ on the $n$th Morava $E$-theory of spaces.
We study a relationship between the actions 
of the Hecke algebras 
$\mathcal{H}_n$ and $\mathcal{H}_{n+1}$
on the $n$th Morava $E$-theory,
and show that the $\mathcal{H}_{n+1}$-module
structure is obtained 
from the $\mathcal{H}_n$-module structure
by the restriction 
along an algebra homomorphism from
$\mathcal{H}_{n+1}$ to $\mathcal{H}_n$. 
\end{abstract}

\maketitle

\section{Introduction}

The stable homotopy category localized 
at a prime number $p$
has a filtration called the chromatic filtration.
Subquotients of the filtration
are said to be monochromatic categories.
It is known that 
the $n$th monochromatic category
is equivalent to the $K(n)$-local category,
where $K(n)$ is the $n$th Morava $K$-theory. 
Thus, 
the stable homotopy category localized at $p$
can be considered to be constructed 
from various $K(n)$-local categories
\cite{Morava, MRW, Ravenel-green, 
DHS, Hopkins-Smith, Ravenel-red}.
Therefore,
it is important to understand the
relationship between $K(n)$-local categories.

The study for relationships between $K(n)$-local categories
of different heights 
is called transchromatic homotopy theory and
it has been done by many authors.
Greenlees-Sadofsky~\cite{Greenlees-Sadofsky}
showed that
the Tate cohomology of $K(n)$ is trivial
for any finite group,
which implies that the Tate cohomology
of a complex-oriented and $v_n$-periodic spectrum
is $v_{n-1}$-periodic.
Hovey-Sadofsky~\cite{Hovey-Sadofsky}
generalized this result and showed that
the Tate cohomology of an $E(n)$-local spectrum
is $E(n-1)$-local,
where $E(n)$ is the $n$th Johnson-Wilson spectrum.
Ando-Morava-Sadofsky~\cite{AMS}
constructed a splitting of a completion
of the $\mathbb{Z}/p$-Tate cohomology of $E(n)$ into a wedge of $E(n-1)$.
In \cite{Torii0}
we generalized this result and constructed a splitting
of a completion of the $(\mathbb{Z}/p)^k$-geometric 
fixed points spectrum of $E_n$ into a wedge of $E_{n-k}$,
where $E_n$ is the $n$th Morava $E$-theory.
These works are related to the chromatic 
splitting conjecture (cf.~\cite{Hovey}).

Hopkins-Kuhn-Ravenel character theory~\cite{HKR1, HKR2}
is another direction of study
for relationships between $K(n)$-local categories.
It describes the $E_n$-cohomology of classifying spaces
of finite groups tensored with the field $\mathbb{Q}$
of rational numbers
in terms of generalized group characters.
This can be considered as a result for 
a relationship between the $K(n)$-local category and 
the $K(0)$-local category. 
In \cite{AMS, Torii1, Torii3}
we constructed a generalization of Chern character
which is a multiplicative natural transformation
of cohomology theories 
from $E_n$ to an extension of $E_{n-1}$.
In \cite{Torii6}
we studied the generalized Chern character
of classifying spaces of finite groups
by using Hopkins-Kuhn-Ravenel character theory.
Stapleton~\cite{Stapleton1,Stapleton2}
generalized Hopkins-Kuhn-Ravenel character theory,
and gave a description of the Borel equivariant
$E_n$-cohomology tensored with an extension
of $L_{K(t)}E_{n}$ for $0\le t<n$,
where $L_{K(t)}$ is the Bousfield localization functor
with respect to $K(t)$,
by using the language of $p$-divisible groups.
This work can be considered as a study
for a relationship between 
the $K(n)$-local category and
the $K(t)$-local category for $0\le t< n$.
Lurie~\cite{Lurie} further generalized this work
and put it in a more general framework
by introducing a notion of
tempered cohomology theories.


The $K(n)$-local category can be studied by 
the $n$th Morava $E$-theory $E_n$
with action of the extended Morava stabilizer
group $\mathbb{G}_n$.
For example, 
$E_n$-based Adams spectral sequences
strongly converge to 
the homotopy groups of $K(n)$-local spectra.
The $E_2$-pages of the spectral sequences
are described in terms of the derived functor ${\rm Ext}$
in the category of twisted continuous 
$\mathbb{G}_n$-modules over $\pi_*(E_n)$
under some conditions.
Thus, 
we would like to understand a relationship
between the $n$th Morava $E$-theory $E_n$
with $\mathbb{G}_n$-action
and the $(n+1)$st Morava $E$-theory $E_{n+1}$
with $\mathbb{G}_{n+1}$-action.

For this purpose,
we constructed a generalization of the classical
Chern character
\[ {\rm ch}: E_{n+1}\longrightarrow \mathbb{B}_n \]
in \cite{Torii1, Torii3}, 
which is a map of ring spectra from the chromatic height $(n+1)$
complex orientable spectrum $E_{n+1}$
to a chromatic height $n$ complex orientable spectrum $\mathbb{B}_n$.
This type of generalization of the classical Chern character
was first considered by Ando-Morava-Sadofsky~\cite{AMS}.
In \cite{Torii5} we showed that
the generalized Chern character ${\rm ch}$
can be lifted to a map of $E_{\infty}$-ring spectra,
and the spectrum $\mathbb{B}_n$ can be identified with
$L_{K(n)}(E_n\wedge E_{n+1})$.
Under this identification,
the generalized Chern character 
${\rm ch}: E_{n+1}\to \mathbb{B}_n$ 
is the inclusion into the second factor of the
smash product in $\mathbb{B}_n$.
We define a map of $E_{\infty}$-ring spectra
\[ {\rm inc}: E_n \longrightarrow \mathbb{B}_n \]
to be the inclusion
into the first factor of the smash product in
$\mathbb{B}_n$.

Hecke operators are defined by using power operations.
In \cite{Goerss-Hopkins} Goerss-Hopkins
showed that Morava $E$-theory supports
an $E_{\infty}$-ring spectrum structure
which is unique up to homotopy.
The $E_{\infty}$-ring structure on $E_n$
gives rise to power operations, and hence Hecke operators.
See, for example, Ando-Hopkins-Strickland~\cite{AHS},
Rezk~\cite{Rezk2}, and Ganter~\cite{Ganter}.
But Hecke operators in Morava $E$-theory
were first defined by Ando in \cite{Ando1, Ando2},
where he used power operations in the complex cobordism
to produce power operations in Morava $E$-theory.
In the end
there is a natural action of a kind of Hecke algebra
$\mathcal{H}_n$ on $E_n^0(X)$ for any space $X$.

In this paper we consider a relationship between
the Hecke operators in $E_n$ and those in $E_{n+1}$. 
For this purpose,
we construct Hecke operators
in $\mathbb{B}_n$-theory
\[ \widetilde{\rm T}_M^{\mathbb{B}}:
   \mathbb{B}_n^0(X)\longrightarrow \mathbb{B}_n^0(X)\]
for any space $X$ and each finite abelian $p$-group $M$
with $p$-rank $\le n+1$.
The operations $\widetilde{\rm T}^{\mathbb{B}}_M$
are natural extensions of the Hecke operators in $E_{n+1}$.
The following is our first main theorem.

\begin{theorem}
[{Theorem~\ref{thm:Hn+1-module-structure-on-BnX}}]
\label{Main-Theorem-1}
Assigning the Hecke operator $\widetilde{\rm T}_M^{\mathbb{B}}$
to a finite abelian $p$-group $M$ 
with $p$-rank $\le n+1$,
there is a natural action of the Hecke algebra
$\mathcal{H}_{n+1}$ on $\mathbb{B}_n^0(X)$
for any space $X$ such that
${\rm ch}: E_{n+1}^0(X)\to \mathbb{B}_n^0(X)$
is a map of $\mathcal{H}_{n+1}$-modules.
\end{theorem}

Let $\mathbb{G}_{n+1}$ be the $(n+1)$st
extended Morava stabilizer group.
In \cite{Torii3} we showed that
there is a natural action of $\mathbb{G}_{n+1}$
on $\mathbb{B}_n^0(X)$ and that
the map ${\rm inc}$ induces
a natural isomorphism
\[ E_n^0(X)\stackrel{\cong}{\to}
   (\mathbb{B}_n^0(X))^{\mathbb{G}_{n+1}} \]
for any spectrum $X$,
where the right hand side is the $\mathbb{G}_{n+1}$-invariant
submodule of $\mathbb{B}_n^0(X)$.
We shall show that
the action of $\mathcal{H}_{n+1}$ on $\mathbb{B}_n^0(X)$
commutes with the action of $\mathbb{G}_{n+1}$.
This implies the following theorem, 
which is our second main theorem.

\begin{theorem}
[{Theorem~\ref{thm:natural-Hn+1-module-strucrure-on-EnX}}]
\label{Main-Theorem-2}
There is a natural action of the Hecke algebra
$\mathcal{H}_{n+1}$ on $E_n^0(X)$
for any space $X$ such that
${\rm inc}: E_n^0(X)\to \mathbb{B}_n^0(X)$
is a map of $\mathcal{H}_{n+1}$-modules.
\end{theorem}

By these two theorems,
we have the $\mathcal{H}_n$-module structure and
the $\mathcal{H}_{n+1}$-module structure on $E_n^0(X)$,
and we shall compare these two structures.
For this purpose, 
we shall construct an algebra homomorphism
\[ \omega: \mathcal{H}_{n+1}\longrightarrow \mathcal{H}_n,\]
and prove the following theorem,
which is our third main theorem.

\begin{theorem}[{Theorem~\ref{thm:comparison-Hn-Hn+1-on-EnX}}]
\label{Main-Theorem-3}
The $\mathcal{H}_{n+1}$-module structure
on $E_n^0(X)$ is obtained from
the $\mathcal{H}_n$-module structure 
by the restriction along $\omega$.
\end{theorem}

In \cite{Torii8}
we studied a relation between
the category of modules 
over the algebra of power operations
in the $n$th Morava $E$-theory 
and that in the $(n+1)$st Morava $E$-theory
by the same technique in this paper.


\begin{overview}\rm
In \S\ref{section:MoravaE}
we first review the construction of Hecke operators in Morava $E$-theory.
We also study a natural action of the extended Morava stabilizer group
$\mathbb{G}_n$ on level structures
on the formal group associated to the Morava $E$-theory $E_n$.
As a result,
we show that the action of the Hecke algebra $\mathcal{H}_n$
on $E_n^0(X)$ commutes with the action of $\mathbb{G}_n$.
In \S\ref{section:Bn-theory}
we construct Hecke operators in $\mathbb{B}_n$-theory.
For this purpose, 
we introduce and study level structures on 
a $p$-divisible group.
We show that functors for level structures of given type
on the $p$-divisible group associated to $\mathbb{B}_n$-theory
are representable by finite products of complete regular local rings.
Using this result,
we construct additive unstable operations in $\mathbb{B}_n$-theory,
which are natural extensions of the power operations 
in $E_{n+1}$-theory.
After that,
we define Hecke operators in $\mathbb{B}_n$-theory
by assembling these operations,
show that the Hecke algebra $\mathcal{H}_{n+1}$
naturally acts on $\mathbb{B}_n^0(X)$ for any space $X$,
and prove Theorem~\ref{Main-Theorem-1}
(Theorem~\ref{thm:Hn+1-module-structure-on-BnX}).
Furthermore, 
we study a natural $\mathbb{G}_{n+1}$-action
on representing rings 
for level structures on the $p$-divisible group
associated to $\mathbb{B}_n$-theory.
We show that
the action of the Hecke algebra $\mathcal{H}_{n+1}$
on $\mathbb{B}_n^0(X)$ 
commutes with the action of the extended
Morava stabilizer group $\mathbb{G}_{n+1}$,
and prove Theorem~\ref{Main-Theorem-2}
(Theorem~\ref{thm:natural-Hn+1-module-strucrure-on-EnX}).
In \S\ref{section:comparison-Hecke-operators}
we compare the Hecke operators
in $E_n$-theory with those in $E_{n+1}$-theory 
via $\mathbb{B}_n$-theory.
We discuss the restriction of
the Hecke operators in $\mathbb{B}_n$-theory
to $E_n$-theory,
and show that it
can be written in terms of the Hecke operators
in $E_n$-theory.
Using this result, 
we construct a ring homomorphism
$\omega$ from $\mathcal{H}_{n+1}$ to $\mathcal{H}_n$,
and prove Theorem~\ref{Main-Theorem-3}
(Theorem~\ref{thm:comparison-Hn-Hn+1-on-EnX}).
\end{overview}



\begin{notation}\rm

In this paper we fix a prime number $p$ and a positive integer $n$.
Let $\mathbb{Z}_p$ be the ring of $p$-adic integers,
and let $\mathbb{Q}_p$ be its fraction field.
For a set $S$, 
we denote by $|S|$ the cardinality of $S$,
and by $\mathbb{Z}[S]$
the free $\mathbb{Z}$-module generated by $S$.
For a finite abelian $p$-group $M$,
we denote the $p$-rank of $M$
by $p$-rank($M$).
We let $M[p^r]$ be the kernel of 
the multiplication map $p^r:M\to M$
for a nonnegative integer $r$.

We write $\mathcal{CL}$  
for the category of complete Noetherian local rings
with residue field of characteristic $p$
and local ring homomorphisms.
For $R\in \mathcal{CL}$, 
we let $\mathcal{CL}_R$ be the under category $\mathcal{CL}_{R/}$.
For a (formal) scheme $\mathbf{X}$ 
over $R\in\mathcal{CL}$ and 
a morphism $f: R\to S$ in $\mathcal{CL}$,
we denote by $\mathbf{X}_S$ or $f^*\mathbf{X}$
the base change of $\mathbf{X}$ along $f$.
We let $\mathbf{G}^0$ be the identity component of 
a $p$-divisible group $\mathbf{G}$ over
$R\in\mathcal{CL}$.
\end{notation}



\section{Hecke operators in Morava $E$-theory}
\label{section:MoravaE}

In \cite{Goerss-Hopkins} Goerss-Hopkins showed that
Morava $E$-theory supports 
an $E_{\infty}$-ring structure which is unique up to homotopy.
This induces power operations in 
Morava $E$-theory of spaces, and 
we can define Hecke operators by using power operations.
In this section we review the construction of Hecke operators
in Morava $E$-theory
and study their compatibility with the action 
of the extended Morava stabilizer group.

\subsection{Morava $E$-theory}
\label{subsection:Morava-E-theory}

In this subsection
we fix notation of Morava $E$-theory.

The $n$th Morava $E$-theory $E_n$
is an even periodic commutative ring spectrum.
The degree $0$ coefficient ring $E_n^0=\pi_0(E_n)$ 
is given by
\[ W(\mathbb{F}_{p^n})\power{w_1,\ldots,w_{n-1}}, \]
where $W(\mathbb{F}_{p^n})$ is
the ring of Witt vectors with
coefficients in the finite field $\mathbb{F}_{p^n}$
with $p^n$ elements.
Since $E_n$ is even periodic,
we have an associated degree $0$
one dimensional formal Lie group 
$\mathbf{F}_n$ over $E_n^0$,
which is a universal deformation 
of the Honda formal group $\overline{\mathbf{F}}_n$ 
over $\mathbb{F}_{p^n}$.
We denote by $\mathbb{G}_n$ the $n$th extended Morava stabilizer group,
which is a semidirect product
of the automorphism group $\mathbb{S}_n$
of the Honda formal group 
$\overline{\mathbf{F}}_n$ over $\mathbb{F}_{p^n}$
with the Galois group ${\rm Gal}(\mathbb{F}_{p^n}/\mathbb{F}_p)$
of $\mathbb{F}_{p^n}$ over the prime field $\mathbb{F}_p$:
\[ \mathbb{G}_n=\,{\rm Gal}(\mathbb{F}_{p^n}/\mathbb{F}_p)
                  \ltimes \mathbb{S}_n.\]


\subsection{Level structures on formal Lie groups}
\label{subsection:level-structure}

In this subsection we review
level structures on one dimensional formal Lie groups
(cf.~\cite{Strickland}).

Let $\mathbf{X}$ be a one dimensional formal Lie group 
over $R\in\mathcal{CL}$ which has finite height.
First, we recall the definition of divisors on $\mathbf{X}$.
A divisor on $\mathbf{X}$ is a closed subscheme
which is finite and flat over $R$.
We can associate to a section $s\in \mathbf{X}(R)$ 
a divisor $[s]$ on $\mathbf{X}$.
If we take a coordinate $x$ of $\mathbf{X}$,
then the divisor $[s]$ is the closed subscheme
${\rm Spf}(R\power{x}/(x-x(s)))$,
where $x(s)$ is the image of $x$
under the map $s^*: R\power{x}\to R$.

Suppose that we have a homomorphism
$\phi: M\to \mathbf{X}(R)$ of abelian groups,
where $M$ is a finite abelian $p$-group.
We can define a divisor $[\phi(M)]$ on $\mathbf{X}$ 
by
\[ [\phi(M)]=\sum_{m\in M}\ [\phi(m)].\]
If we take a coordinate $x$ of $\mathbf{X}$,
then the divisor $[\phi(M)]$
is the closed subscheme ${\rm Spf}(R\power{x}/(f(x)))$,
where $f(x)=\prod_{m\in M}(x-x(\phi(m)))$.



Next, we recall the definition
of level structures on $\mathbf{X}$. 
For a nonnegative integer $r$,
we have a closed subgroup scheme 
$\mathbf{X}[p^r]$
which is the kernel of 
the map $p^r: \mathbf{X}\to \mathbf{X}$.
A homomorphism $\phi: M\to \mathbf{X}(R)$ 
is said to be a level $M$-structure
if the divisor $[\phi(M[p])]$
is a closed subscheme of $\mathbf{X}[p]$.
Note that 
$[\phi(M[p^r])]$ is a closed subgroup scheme of 
$\mathbf{X}[p^r]$ for all $r\ge 0$
if $\phi: M\to \mathbf{X}(R)$ 
is a level $M$-structure
by \cite[Proposition~32]{Strickland}.

\begin{definition}
\label{def:level-structure-formalgroup}
\rm 
We define a functor 
\[ {\rm Level}(M,\mathbf{X}) \]
from $\mathcal{CL}_R$ to the category of sets
by assigning to $S\in\mathcal{CL}_R$
the set of all level $M$-structures on $\mathbf{X}_S$.
\end{definition}

By \cite[Proposition~22]{Strickland},
the functor ${\rm Level}(M,\mathbf{X})$
is representable,
that is,
there exists $D(M,\mathbf{X})\in\mathcal{CL}_R$ such that
\[ {\rm Level}(M,\mathbf{X})={\rm Spf}(D(M,\mathbf{X})).\]
Note that 
$D(M,\mathbf{X})$ is a finitely generated $R$-module.
Furthermore,
$D(M,\mathbf{X})$ is a regular local ring 
if $\mathbf{X}$ is a universal deformation
of the formal group on the closed point
(see \cite[Theorem~23]{Strickland}).

Next, we consider functoriality of $D(M,\mathbf{X})$
with respect to $M$.
If $\phi:M\to \mathbf{X}(R)$ is a level $M$-structure
and $u: N\to M$ is a monomorphism of abelian $p$-groups,
then the composition $\phi\circ u: N\to \mathbf{X}(R)$
is a level $N$-structure.
Hence we obtain a natural transformation
${\rm Level}(M,\mathbf{X})\longrightarrow
   {\rm Level}(N,\mathbf{X})$,
which induces a map 
\[ D(N,\mathbf{X}) \longrightarrow D(M,\mathbf{X})\]
in $\mathcal{CL}_R$.

Finally,
we recall quotient formal Lie groups.
Let $\phi:M\to \mathbf{X}(R)$ be 
a level $M$-structure on $\mathbf{X}$.
Then the divisor $[\phi(M)]$ is a finite subgroup scheme
of $\mathbf{X}$,
and we can construct a quotient formal Lie group 
$\mathbf{X}/[\phi(M)]$ over $R$
(see \cite[Theorem~19]{Strickland}).



\subsection{The action of $\mathbb{G}_n$ on $D(M,\mathbf{F}_n)$}
\label{subsection:level-structure-Morava-E}

In this subsection we study 
an action of the extended Morava stabilizer group
$\mathbb{G}_n$ on 
the representing ring $D(M,\mathbf{F}_n)$
of level structures on $\mathbf{F}_n$.

First, we recall that 
the action of $\mathbb{G}_n$ on 
the Honda formal group $\overline{\mathbf{F}}_n$ 
over $\mathbb{F}_{p^n}$ extends to
an action on its universal deformation
$\mathbf{F}_n$ over $E_n^0$. 
We will show that 
this action further extends to an action on 
the functor ${\rm Level}(M,\mathbf{F}_n)$ and 
its representing ring $D(M,\mathbf{F}_n)$
for any abelian $p$-group $M$
(cf.~\cite[\S7]{Torii6}).

Let $\phi: M\to \mathbf{F}_n(R)$
be a level $M$-structure on $\mathbf{F}_n$
over $R\in\mathcal{CL}_{E_n^0}$.
For any $g\in \mathbb{G}_n$,
the composite 
\[ M\stackrel{\phi}{\longrightarrow} 
   \mathbf{F}_n(R)\stackrel{g}{\longrightarrow}
   \mathbf{F}_n(R') \]
is a level $M$-structure on $\mathbf{F}_n$ over 
$R'\in \mathcal{CL}_{E_n^0}$,
where $R'$ is a local $E_n^0$-algebra 
given by the composite $E_n^0\stackrel{g}{\to} E_n^0\to R$.
This induces a natural transformation
\[   g: {\rm Level}(M,\mathbf{F}_n) 
        \longrightarrow
     {\rm Level}(M,\mathbf{F}_n) \]
over $g: {\rm Spf}(E_n^0)\to {\rm Spf}(E_n^0)$.
Hence we obtain an action of $\mathbb{G}_n$ on 
the functor ${\rm Level}(M,\mathbf{F}_n)$
which is compatible with the action on ${\rm Spf}(E_n^0)$.
In other words,
we obtain an action of $\mathbb{G}_n$
on its representing ring $D(M,\mathbf{F}_n)$
which is compatible with the action on $E_n^0$.

By the construction of the action of $\mathbb{G}_n$
on level structures,
it is easy to obtain the following lemma.

\begin{lemma}\label{lemma:restriction-g-level-str-commutativity}
For any monomorphism $N\to M$ of finite abelian $p$-groups,
the induced ring homomorphism
$D(N,\mathbf{F}_n)\to D(M,\mathbf{F}_n)$
is $\mathbb{G}_n$-equivariant.
\end{lemma}



\subsection{Power operations in Morava $E$-theory}
\label{subsection:power_operation_MoravaE}

In this subsection
we recall the construction
of power operations in Morava $E$-theory
and show that they are equivariant
with respect to the action of 
the extended Morava stabilizer group.

For a space $X$,
we denote by $\Sigma^{\infty}_+X$
the suspension spectrum of $X$
with a disjoint base point.
By \cite{Goerss-Hopkins},
Morava $E$-theory $E_n$
is an $E_{\infty}$-ring spectrum 
which is unique up to homotopy.
In particular, we have a structure map
\[ \nu_k: \Sigma^\infty_+E\Sigma_k
    \subrel{\Sigma_k}{\wedge}
   E_n^{\wedge k} \longrightarrow E_n,\]
where $\Sigma_k$ is the symmetric group of degree $k$,
and $E\Sigma_k$ is its universal space.
For a map
$f: \Sigma^{\infty}_+X\to E_n$,
by the composition
\[ \Sigma^{\infty}_+(E\Sigma_k\subrel{\Sigma_k}{\times}
   X^{\times k})\cong
   \Sigma^{\infty}_+E\Sigma_k\subrel{\Sigma_k}{\wedge}
   \Sigma_+X^{\wedge k}
   \stackrel{1\wedge f^{\wedge k}}{\hbox to 13mm{\rightarrowfill}}
   \Sigma^{\infty}_+E\Sigma_k\subrel{\Sigma_k}{\wedge}
   E_n^{\wedge k}
   \stackrel{\nu_{k}}{\hbox to 10mm{\rightarrowfill}}E_n,\]
we obtain a natural transformation
\[ {\rm p}_k:E_n^0(X)\longrightarrow 
   E_n^0(E\Sigma_k\subrel{\Sigma_k}{\times}X^{\times k}).\]
The inclusion map of the diagonal 
$\Delta: X\to X^{\times k}$ induces
a map
$(1\times\Delta)^*:
   E_n^0(E\Sigma_k\subrel{\Sigma_k}{\times}X^{\times k})
   \to
   E_n^0(E\Sigma_k\subrel{\Sigma_k}{\times}X)\cong
   E_n^0(B\Sigma_k)\subrel{E_n^0}{\otimes} E_n^0(X)$.
A power operation 
\[ {\rm P}_{k}:  E_n^0(X)\longrightarrow 
   E_n^0(B\Sigma_k)\subrel{E_n^0}{\otimes}
   E_n^0(X) \]
is defined by
\[ {\rm P}_{k}= (1\times\Delta)^*\circ {\rm p}_k.\]

Let $G$ be a finite group of order $k$.
We fix a bijection between $G$
and the set $\{1,2,\ldots,k\}$. 
The left multiplication by an element in $G$
induces an automorphism of $\{1,2,\ldots,k\}$.
Hence we can regard $G$ as a subgroup of $\Sigma_k$.
The inclusion $i:G\hookrightarrow \Sigma_k$
induces a map $Bi: BG\to B\Sigma_k$
of classifying spaces.
Note that the homotopy class of $Bi$
is independent of the choice 
of a bijection between $G$ and $\{1,2,\ldots,k\}$.
We define a power operation 
\[ {\rm P}_{G}:  E_n^0(X)\longrightarrow 
   E_n^0(BG)\subrel{E_n^0}{\otimes} E_n^0(X) \]
by the composition
\[ {\rm P}_G= (Bi^*\otimes 1)\circ {\rm P}_k.\]

For a finite abelian $p$-group $M$,
we set $M^*={\rm Hom}(M,S^1)$,
where $S^1$ is the circle group.
By \cite[Proposition~2.9]{Greenlees-Strickland},
the functor ${\rm Hom}(M,\mathbf{F}_n)$
over ${\rm Spf}(E_n^0)$
is represented by $E_n^0(BM^*)$,
that is, 
\[ {\rm Hom}(M,\mathbf{F}_n)=\mbox{\rm Spf}(E_n^0(BM^*)).\]
Since the closed subscheme ${\rm Level}(M,\mathbf{F}_n)$
of ${\rm Hom}(M,\mathbf{F}_n)$
is represented by $D(M,\mathbf{F}_n)$,
we have a local $E_n^0$-algebra homomorphism
\[ {\rm R}: E_n^0(BM^*)\longrightarrow D(M,\mathbf{F}_n).\]  

Using the $E_n^0$-algebra $D(M,\mathbf{F}_n)$,
we define an extension of Morava $E$-theory.

\begin{definition}\rm
For a space (or a spectrum) $X$, we set 
\[ E(M)_n^*(X)= D(M,\mathbf{F}_n)\subrel{E_n^0}{\otimes}E_n^*(X).\]
Since $D(M,\mathbf{F}_n)$ 
is a finitely generated free $E_n^0$-module,
the functor $E(M)_n^*(X)$ is a generalized cohomology theory.
\end{definition}

By the diagonal action,
the group $\mathbb{G}_n$ acts on $E(M)_n^*(X)$ as multiplicative
stable cohomology operations.
For a monomorphism $N\to M$ of finite abelian $p$-groups,
we have a multiplicative stable cohomology operation
$E(N)_n^*(X)\to E(M)_n^*(X)$,
which is $\mathbb{G}_n$-equivariant 
by Lemma~\ref{lemma:restriction-g-level-str-commutativity}.

\begin{definition}\rm
We define an operation
\[ \Psi_M^{E_n}: 
   E_n^0(X)\longrightarrow E(M)_n^0(X)\]
by the composition
\[ \Psi_M^{E_n}= ({\rm R}\otimes 1)\circ {\rm P}_{M^*}.\]
The operation
$\Psi_M^{E_n}$ is a ring operation (see \cite[Theorem~3.4.4]{Ando1}
or \cite[Lemma~3.10]{AHS}).
\end{definition}

When $X$ is a one point space, 
we obtain a ring homomorphism
\[ \Psi_M^{E_n}: E_n^0\longrightarrow E(M)_n^0=D(M,\mathbf{F}_n).\]
This map classifies the quotient formal group  
$i^*\mathbf{F}_n/[\phi_{\rm univ}(M)]$,
where $i: E_n^0\to D(M,\mathbf{F}_n)$ is
the structure map and $\phi_{\rm univ}$ is the universal
level $M$-structure on $i^*\mathbf{F}_n$.
When $X$ is the infinite complex projective space 
$\mathbb{CP}^{\infty}$,
the map
$\Psi_M^{E_n}: E_n^0(\mathbb{CP}^{\infty})
\to E(M)_n^0(\mathbb{CP}^{\infty})$ induces 
an isogeny
\[ \Psi_M^{E_n}:
   i^*\mathbf{F}_n\longrightarrow \Psi_M^*\mathbf{F}_n\]
with 
$\ker\Psi_M^{E_n}=[\phi_{\rm univ}(M)]$
(see \cite[\S\S3.3 and 3.4]{Ando1}
or \cite[Proposition~3.21]{AHS}).

Next, we consider compatibility of 
power operations $\Psi_M^{E_n}$
with the action of the extended Morava stabilizer group
$\mathbb{G}_n$.

\begin{proposition}\label{prop:g-Psi-commutation-En}
The operation
$\Psi_M^{E_n}: E_n^0(X)\longrightarrow E(M)_n^0(X)$
is $\mathbb{G}_n$-equivariant.
\end{proposition}

\proof
It suffices to show that
$P_{M^*}$ and $R$ are $\mathbb{G}_n$-equivariant.
Recall that $\mathbb{G}_n$ acts on $E_n$ as $E_{\infty}$-ring spectrum maps.
This implies that $P_{M^*}$ is $\mathbb{G}_n$-equivariant.
By the construction of the action of $\mathbb{G}_n$,
the map of functors
${\rm Level}(M,\mathbf{F}_n)\to {\rm Hom}(M,\mathbf{F}_n)$
is $\mathbb{G}_n$-equivariant.
Thus, 
$R: E_n^0(BM^*)\to D(M,\mathbf{F}_n)$ is also $\mathbb{G}_n$-equivariant.
\qqq

\subsection{Hecke operators in Morava $E$-theory}

In this subsection we review Hecke operators
in Morava $E$-theory
(see, for example, \cite{Ganter,Rezk,Ando2,AHS}).
We also study their compatibility 
with the action of the extended Morava stabilizer
group.

First, 
we recall the construction of Hecke operators 
in Morava $E$-theory.
Let
$\Lambda^n=(\mathbb{Q}_p/\mathbb{Z}_p)^n$.
We set 
\[ D_n(r)=D(\Lambda^n[p^r],\mathbf{F}_n). \]
The inclusion 
$\Lambda^n[p^r]\to \Lambda^n[p^{r+1}]$
of finite abelian $p$-groups
induces a local $E_n^0$-algebra homomorphism
$D_n(r)\to D_n(r+1)$.
We define an $E_n^0$-algebra $D_n$ to be
a colimit of $\{D_n(r)\}_{r\ge 0}$:
\[ D_n=\ \subrel{r}{\rm colim} D_n(r).\]

For a subgroup $M$ of $\Lambda^n$,
we let ${\rm I}_M: D(M,\mathbf{F}_n)\to D_n$
be the $E_n^0$-algebra homomorphism
induced by the inclusion.
We define an operation 
\[ \Phi_M^{E_n}: E_n^0(X)\longrightarrow 
   D_n\subrel{E_n^0}{\otimes}E_n^0(X)\]
by the composition 
\[ \Phi_M^{E_n}= ({\rm I}_M\otimes 1)\circ \Psi_M^{E_n}.\]

We would like to have operations
from $E_n^0(X)$ to $E_n^0(X)$.
For this purpose,
we consider a sum of operations $\Phi_M^{E_n}$ for 
all subgroups of $\Lambda^n$ isomorphic to 
a given finite abelian $p$-group $M$.
Then we can show that it factors through $E_n^0(X)$
by an invariance with respect
to the Galois group of $D_n$ over $E_n^0$.

For a finite abelian $p$-group $M$,
we let $m_n(M)$ be the set of all subgroups
of $\Lambda^n$ isomorphic to $M$:
\[ m_n(M)=\{M'\le \Lambda^n|\ M'\cong M\}.\]
Note that $m_n(M)$ is a finite set
and that 
$m_n(M)=\emptyset$
if $p$-rank$(M)>n$.
By \cite[\S14.2]{Rezk},
there is an unstable additive cohomology operation
\[ \widetilde{\rm T}_M: E_n^0(X)\longrightarrow E_n^0(X) \]
such that
\[ (i\otimes 1)\circ \widetilde{\rm T}_M=
   \sum_{M'\in m_n(M)} \Phi_{M'}^{E_n}, \]
where $i:E_n^0\to D_n$ is the $E_n^0$-algebra
structure map of $D_n$.

In the same way
there is an unstable additive operation
\[ \widetilde{\rm T}\left({p^r}\right): 
   E_n^0(X)\longrightarrow E_n^0(X)\]
for a nonnegative integer $r$ such that
\[ (i\otimes 1)\circ\widetilde{\rm T}\left({p^r}\right)=
   \sum_{\mbox{\scriptsize$
         \begin{array}{c}
         M' \le \Lambda^n\\          
         |M'|=p^r\\
         \end{array}$}} \Phi_{M'}^{E_n},\]
where the sum ranges
over the subgroups $M'$ of $\Lambda^n$
such that the order $|M'|$ is $p^r$.

Now, we recall an algebra of Hecke operators
associated to $(\Delta_n,\Gamma_n)$,
where 
$\Delta_n={\rm End}(\Lambda^n)\cap {\rm GL}_n(\mathbb{Q}_p)$
and $\Gamma_n={\rm GL}_n(\mathbb{Z}_p)$
(cf.~\cite[Chapter~3]{Shimura} and
\cite[\S14.1]{Rezk}).
We denote by $\mathbb{Z}[\Delta_n]$
the monoid ring of $\Delta_n$ over $\mathbb{Z}$,
and by $\mathbb{Z}[\Delta_n/\Gamma_n]$
the left $\mathbb{Z}[\Delta_n]$-module 
spanned by cosets.
We define an algebra $\mathcal{H}_n$ 
to be the endomorphism ring of the
left $\mathbb{Z}[\Delta_n]$-module $\mathbb{Z}[\Delta_n/\Gamma_n]$:
\[ \mathcal{H}_n={\rm End}_{\mathbb{Z}[\Delta_n]}
   (\mathbb{Z}[\Delta_n/\Gamma_n]).\]

Note that $\mathcal{H}_n$
is isomorphic to $\mathbb{Z}[\Gamma_n\backslash\Delta_n/\Gamma_n]$
as $\mathbb{Z}$-modules,
where $\mathbb{Z}[\Gamma_n\backslash\Delta_n/\Gamma_n]$
is the free $\mathbb{Z}$-module generated by double cosets.
We identify $\Gamma_n\backslash\Delta_n/\Gamma_n$
with the set of
all isomorphism classes of finite abelian $p$-groups
of $p$-rank $\le n$
by associating $\ker\alpha$
to the double coset $\Gamma_n\alpha\Gamma_n$.
For a finite abelian $p$-group $M$ 
of $p$-rank $\le n$,
the associated endomorphism
$\widetilde{M}: \mathbb{Z}[\Delta_n/\Gamma_n]\to
\mathbb{Z}[\Delta_n/\Gamma_n]$
is given by
$\Gamma_n\mapsto \sum\alpha\Gamma_n$,
where the sum is taken over
$\alpha\Gamma_n\in\Delta_n/\Gamma_n$
such that $\ker\alpha\cong M$.

The composition in $\mathcal{H}_n$ is given as follows.
For finite abelian $p$-groups $L,M,N$ with $p$-rank $\le n$,
we let $K_n(M,N;L)$ be the set of all
the pairs $(M',L')$ of subgroups of $\Lambda^n$
such that
$M'\subset L'$, $L'\cong L$,
$M'\cong M$, and $L'/M'\cong N$.
We set 
\[ c_n(M,N;L)=\frac{|K_n(M,N;L)|}{|m_n(L)|}.\]
The composition in $\mathcal{H}_n$ is given by
\[ \widetilde{M}\cdot \widetilde{N}= 
   \sum\ c_n(M,N;L)\ \widetilde{L} \]
(cf.~\cite[Proposition~3.15]{Shimura}).
By \cite[Theorem~3.20]{Shimura},
$\mathcal{H}_n$ is a polynomial ring over $\mathbb{Z}$
generated by the endomorphisms
associated to 
$(\mathbb{Z}/p\mathbb{Z})^{\oplus k}$ for $1\le k\le n$.
In particular,
$\mathcal{H}_n$ is a commutative algebra over $\mathbb{Z}$.

In \cite[\S14]{Rezk}
it is shown that
there is a natural action of $\mathcal{H}_n$ on $E_n^0(X)$  
for any space $X$.
Assigning to $\widetilde{M}$ the Hecke operator
$\widetilde{\rm T}_M$,
we obtain a natural $\mathcal{H}_n$-module structure
on $E_n^0(X)$ for any space $X$
by \cite[Proposition~14.3]{Rezk}.

Now,
we consider compatibility of
Hecke operators $\widetilde{\rm T}_M$ and
$\widetilde{\rm T}\left(p^r\right)$
with the action of the extended
Morava stabilizer group $\mathbb{G}_n$.

\begin{proposition}
\label{prop:commutation-Hecke-stabilizer-En}
The operations
$\widetilde{\rm T}_M$ and 
$\widetilde{\rm T}\left({p^r}\right)$
are $\mathbb{G}_n$-equivariant.
\end{proposition}

\proof
%
The proposition follows from 
Lemma~\ref{lemma:restriction-g-level-str-commutativity}
and
Proposition~\ref{prop:g-Psi-commutation-En}.
\qqq

 

\begin{corollary}
The action of $\mathcal{H}_n$ on $E_n^0(X)$
commutes with the action of $\mathbb{G}_n$.
\end{corollary}

\section{Hecke operators in $\mathbb{B}_n$-theory}
\label{section:Bn-theory}

In \cite{Torii3} we defined an even periodic ring spectrum
$\mathbb{B}_n$ that connects the Morava
$E$-theories $E_n$ and $E_{n+1}$. 
In this section
we construct unstable additive operations 
in $\mathbb{B}_n$,
which are extensions of power operations
in $E_{n+1}$.
Using these operations,
we define Hecke operators in $\mathbb{B}_n$-theory.

\subsection{$\mathbb{B}_n$-theory}

In this subsection 
we recall an even periodic
commutative ring spectrum $\mathbb{B}_n$ 
which is an extension of $E_n$ and $E_{n+1}$.
We show that the degree $0$ coefficient ring $\mathbb{B}_n^0$ 
classifies isomorphisms between 
$p$-divisible groups associated to $E_n$ and $E_{n+1}$.

First, we introduce a spectrum $\mathbb{A}_n$
which is defined by
\[ \mathbb{A}_n= L_{K(n)}E_{n+1},\]
where $K(n)$ is the $n$th Morava $K$-theory spectrum
and $L_{K(n)}$ is the Bousfield localization functor
with respect to $K(n)$.
The spectrum $\mathbb{A}_n$ is an even periodic commutative ring
spectrum.
In order to distinguish generators of $E_n^0$ and $E_{n+1}^0$,
we write $E_{n+1}^0=
W(\mathbb{F}_{p^{n+1}})\power{u_1,\ldots,u_n}$.
Then the degree $0$ coefficient ring 
of $\mathbb{A}_n$ is given by
\[ \mathbb{A}_n^0=(W((u_n)))_p^{\wedge}
   \power{u_1,\ldots,u_{n-1}},\]
where $W=W(\mathbb{F}_{p^{n+1}})$.

We describe a relationship between
the formal groups of $E_{n+1}$ and $\mathbb{A}_n$.
Let $\mathbf{G}$ be
the $p$-divisible group $\mathbf{F}_{n+1}[p^{\infty}]$
associated to the formal group $\mathbf{F}_{n+1}$
over $E_{n+1}^0$.
We denote by
$\mathbf{G}_{\mathbb{A}}$
the base change of $\mathbf{G}$
along the map 
$E_{n+1}^0\to \mathbb{A}_n^0$.
There is an exact sequence
$0\to \mathbf{G}_{\mathbb{A}}^0\to\mathbf{G}_{\mathbb{A}}\to
\mathbf{G}_{\mathbb{A}}^{\rm et}\to 0$
of $p$-divisible groups over $\mathbb{A}_n^0$,
where $\mathbf{G}_{\mathbb{A}}^0$
is connected and 
$\mathbf{G}_{\mathbb{A}}^{\rm et}$ is \'{e}tale.
The $p$-divisible group associated to 
the formal group of 
$\mathbb{A}_n$
is equivalent to 
$\mathbf{G}_{\mathbb{A}}^0$.

Next,
we introduce a spectrum $\mathbb{B}_n$ 
which is an amalgamation of $E_n$ and $E_{n+1}$.
We define a spectrum $\mathbb{B}_n$ by
\[ \mathbb{B}_n=L_{K(n)}(E_n{\wedge} E_{n+1}).\]

The spectrum $\mathbb{B}_n$ is also an even periodic
commutative ring spectrum. 
We have ring spectrum maps
\[ {\rm inc}: E_n\longrightarrow \mathbb{B}_n,
   \qquad
   {\rm ch}: E_{n+1}\longrightarrow \mathbb{B}_n,\]
where ${\rm inc}$ is the inclusion 
into the first smash factor,
and ${\rm ch}$ is the inclusion 
into the second smash factor of $\mathbb{B}_n$.
Since $\mathbb{B}_n$ is $K(n)$-local,
the map ${\rm ch}$ induces a map of ring spectra
\[ {\rm j}: \mathbb{A}_n\longrightarrow \mathbb{B}_n.\]
In the following of this paper
we will use same symbols for induced ring homomorphisms
on $\pi_0$.

We regard 
the degree $0$ coefficient ring $\mathbb{B}_n^0$
as an extension of $\mathbb{A}_n^0$
along ${\rm j}$.
Then $\mathbb{B}_n^0$
is a complete regular local ring
with maximal ideal $(p,u_1,\ldots,u_{n-1})$
and that its residue field is a Galois
extension of $\mathbb{F}_{p^{n+1}}((u_n))$
with Galois group isomorphic to 
the extended Morava stabilizer group 
$\mathbb{G}_n$
(see \cite[\S4]{Torii3}).

Now, we study a functor represented 
by the degree $0$ coefficient ring $\mathbb{B}_n^0$.
Let $\mathbf{H}$ be 
the $p$-divisible group $\mathbf{F}_n[p^{\infty}]$
associated to the formal group $\mathbf{F}_n$
over $E_n^0$.
We suppose that there are two maps
$f: \mathbb{A}_n^0\to R$ and $g: E_n^0\to R$
in $\mathcal{CL}$.
Then we have two connected $p$-divisible groups
$f^*\mathbf{G}_{\mathbb{A}}^0$ and
$g^*\mathbf{H}$ over $R$.  
We denote by
\[ {\rm Iso}(f^*\mathbf{G}_{\mathbb{A}}^0, 
   g^*\mathbf{H})\]
the set of all isomorphisms
of $p$-divisible groups between 
$f^*\mathbf{G}_{\mathbb{A}}^0$ and
$g^*\mathbf{H}$,
and by 
\[ {\rm Hom}_{f,g}^c(\mathbb{B}_n^0,R) \]
the set of all maps of local rings
$h: \mathbb{B}_n^0 \to R$ such that 
$h\circ {\rm j}=f$ and 
$h\circ {\rm inc}=g$.

\begin{lemma}\label{lemma:MUPMUP-variation}
There is a bijection 
between 
${\rm Iso}(f^*\mathbf{G}_{\mathbb{A}}^0, 
g^*\mathbf{H})$
and ${\rm Hom}_{f,g}^c(\mathbb{B}_n^0,R)$.
\end{lemma}

\proof
Let ${\rm MP}$ be the periodic complex bordism spectrum.
We have
$\pi_0(E_n\wedge \mathbb{A}_n)=
E_{n}^0\otimes_{{\rm MP}_0}{\rm MP}_0({\rm MP})
\otimes_{{\rm MP}_0}\mathbb{A}_n^0$.
Since $\mathbb{B}_n\simeq L_{K(n)}(E_n\wedge \mathbb{A}_n)$,
the coefficient ring
$\mathbb{B}_n^0$ is a completion
of $\pi_0(E_n\wedge \mathbb{A}_n)$ 
at the ideal $I_n=(p,u_1,\ldots,u_{n-1})$.
The lemma follows from the fact that
${\rm MP}_0({\rm MP})$ classifies isomorphisms of formal group laws.
\qqq



\subsection{Level structures on $p$-divisible groups}

In this subsection we recall the definition
of level structures on $p$-divisible groups.
If a $p$-divisible group is obtained from
a one dimensional formal group,
then we study a relationship between level structures on them.

Let $\mathbf{X}$ be 
a $p$-divisible group over 
a commutative ring $R$.
We assume that 
the finite subgroup scheme $\mathbf{X}[p^r]$ 
is embeddable in a smooth curve 
in the sense of \cite[(1.2.1)]{KM}
for any $r\ge 0$.

\begin{definition}\rm
Let $M$ be a finite abelian $p$-group.
A homomorphism $\phi: M\to \mathbf{X}(R)$
is said to be a level $M$-structure on $\mathbf{X}$
if the induced homomorphism
$\phi[p^r]:M[p^r]\to \mathbf{X}[p^r](R)$
is an $M[p^r]$-structure on $\mathbf{X}[p^r]$ 
in the sense of \cite[Remark~1.10.10]{KM}
for any $r\ge 0$.
\end{definition}

\begin{definition}\rm
We define a functor 
\[ {\rm Level}(M,\mathbf{X}) \]
over ${\rm Spec}(R)$ by assigning
to an $R$-algebra $U$ 
the set of all
the level $M$-structures on $\mathbf{X}_U$.
By \cite[Lemma~1.3.4]{KM},
the functor ${\rm Level}(M,\mathbf{X})$
is representable.
We denote by $D(M,\mathbf{X})$
the representing $R$-algebra so that
\[ {\rm Level}(M,\mathbf{X})
   ={\rm Spec}(D(M,\mathbf{X})).\]
\end{definition}

Now, we suppose that there exists 
a one dimensional formal Lie group $\mathbf{X}'$
of finite height over $R'\in\mathcal{CL}$
and that $\mathbf{X}$ is obtained from
the associated $p$-divisible group $\mathbf{X}'[p^{\infty}]$
by base change along a ring homomorphism $R'\to R$.
By taking a coordinate $x$ of $\mathbf{X}'$,
we have a $p^r$-series $[p^r](x)\in R'\power{x}$ of 
the associated formal group law.
The Weierstrass preparation theorem
implies that there is a unique monic polynomial 
$g_r(x)\in R'[x]$ 
such that $[p^r](x)$ is a unit multiple
of $g_r(x)$ in $R'\power{x}$.
In this case the subgroup scheme $\mathbf{X}[p^r]$
is given by
\[ \mathbf{X}[p^r]\cong {\rm Spec}(R[x]/(g_r(x))).\] 
Note that 
we can embed $\mathbf{X}[p^r]$ in the smooth curve
${\rm Spec}(R[x])$ for any $r\ge 0$.

We shall describe the above definition
of level $M$-structure
more explicitly.
For a section $s\in \mathbf{X}[p^r](U)$,
we denote by $x(s)$
the image of $x$ under the 
$R$-algebra homomorphism $s^*: R[x]/(g_r(x))\to U$.
By \cite[Lemma~1.10.11]{KM},
there is at most one closed subgroup scheme $K$  
of $\mathbf{X}[p^r]_U$ such that $(K,\phi)$
is an $M[p^r]$-structure on $\mathbf{X}[p^r]_U$.
If such a $K$ exists, 
then the set $\{\phi(m)|\ m\in M[p^r]\}$ is a 
full set of sections of $K$ by \cite[Theorem~1.10.1]{KM},
and
hence the product
\begin{align}\label{align:product-divisor} 
\prod_{m\in M[p^r]}(x-x(\phi(m))) 
\end{align}
divides $g_r(x)$ in $U[x]$.
Conversely,
if the product~(\ref{align:product-divisor}) 
divides $g_r(x)$,
then the divisor
\[ K=\sum_{m\in M[p^r]}[\phi(m)] \]
formed in the smooth curve ${\rm Spec}(U[x])$ 
is a subgroup scheme of $\mathbf{X}[p^r]_U$
by using \cite[Proposition~32]{Strickland},
and hence $(K,\phi)$ is an $M[p^r]$-structure
on $\mathbf{X}[p^r]_U$.
Therefore, $\phi$ is a level $M$-structure
if and only if the product~(\ref{align:product-divisor}) 
divides $g_r(x)$ in $U[x]$ for any $r\ge 0$.

Recall that the functor ${\rm Level}(M,\mathbf{X}')$
over ${\rm Spf}(R')$
is represented by 
$D(M,\mathbf{X}')$. 
By the above argument,
we obtain the following proposition.

\begin{proposition}\label{prop:level-representability-specR}
Let $\mathbf{X}'$ be a one dimensional formal Lie group
over $R'\in\mathcal{CL}$ of finite height.
If 
the $p$-divisible group $\mathbf{X}$
is obtained from the associated $p$-divisible group
$\mathbf{X}'[p^{\infty}]$ by base change along
a map $R'\to R$, 
then 
\[ D(M,\mathbf{X})\cong R\otimes_{R'}D(M,\mathbf{X}').\]
\end{proposition}


\subsection{Extensions of $p$-divisible groups}

Let $\mathbf{X}$ be 
a $p$-divisible group over a commutative ring $R$
such that $\mathbf{X}[p^r]$ 
is embeddable in a smooth curve for any $r\ge 0$.
We suppose that there is an exact sequence 
\[ 0\longrightarrow \mathbf{Y}\longrightarrow
                    \mathbf{X}\stackrel{q}{\longrightarrow}
                    \mathbf{E}\longrightarrow 0\] 
of $p$-divisible groups,
where $\mathbf{E}$ is \'{e}tale.
In this subsection we study a relationship
between level structures on $\mathbf{X}$ and $\mathbf{Y}$.

The map $q: \mathbf{X}\to\mathbf{E}$ induces a homomorphism
$q: \mathbf{X}(R)\to \mathbf{E}(R)$
of abelian groups. 
For a homomorphism
$\phi: M\to \mathbf{X}(R)$,
we set $N=\ker (q\circ \phi)$.
The restriction of $\phi$ to $N$
induces a homomorphism
$\phi': N\to \mathbf{Y}(R)$.

\begin{proposition}
\label{prop:general-level-restriction-height-change}
The homomorphism
$\phi$ is a level $M$-structure on $\mathbf{X}$
if and only if
the restriction $\phi'$ is a level $N$-structure on 
$\mathbf{Y}$.
\end{proposition}

\proof
For any $r\ge 0$, we have a commutative diagram
\[ \begin{array}{ccccccccc}
    0&\longrightarrow& N[p^r]&
      \longrightarrow& M[p^r]&
      \longrightarrow& M[p^r]/N[p^r]&
      \longrightarrow& 0   \\[1mm]
     &&\phantom{\mbox{$\scriptstyle \phi'$}}
       \bigg\downarrow
       \mbox{$\scriptstyle \phi'$}  &
     &\phantom{\mbox{$\scriptstyle \phi$}}\bigg\downarrow
      \mbox{$\scriptstyle \phi$} &
     &\bigg\downarrow& \\[3mm]
    0&\longrightarrow& \mathbf{Y}[p^r](R)&
      \longrightarrow& \mathbf{X}[p^r](R)&
      \longrightarrow& \mathbf{E}[p^r](R),&
                     & \\      
   \end{array} \]
where each row is exact and the right vertical arrow is injective.
The proposition follows in the same way
as in the proof of \cite[Proposition~1.11.2]{KM}.
\qqq


For a level $M$-structure
$\phi: M\to \mathbf{X}(R)$,
we have a finite flat subgroup scheme $[\phi(M)]$
in $\mathbf{X}$.
We denote by $\mathbf{X}/[\phi(M)]$ 
the quotient as a fppf sheaf of abelian groups.
By \cite[Proposition~2.7]{RZ},
$\mathbf{X}/[\phi(M)]$ 
is a $p$-divisible group.
Using the snake lemma,
we obtain the following proposition.

\begin{proposition}
\label{prop:exact-seq-quotient-pdivisible}
Let $\phi: M\to \mathbf{X}(R)$
be a level $M$-structure on $\mathbf{X}$
and let $\phi': N\to \mathbf{Y}(R)$
be the induced level $N$-structure
on $\mathbf{Y}$
by Proposition~\ref{prop:general-level-restriction-height-change},
where $N=\ker (q\circ \phi)$.
Then there exists an exact sequence
\[   0  \longrightarrow  \mathbf{Y}/[\phi'(N)]  
        \longrightarrow  \mathbf{X}/[\phi(M)] 
        \longrightarrow  \mathbf{E}/(M/N)_R 
        \longrightarrow  0\]
of $p$-divisible groups over $R$,
where $(M/N)_R$ is a constant 
group scheme.
\end{proposition}

Next, we consider a decomposition
of the representing ring $D(M,\mathbf{X})$
of the functor ${\rm Level}(M,\mathbf{X})$.
The map $q: \mathbf{X}\to \mathbf{E}$ 
induces a natural transformation
${\rm Level}(M,\mathbf{X})
   \to{\rm Hom}(M,\mathbf{E})$
of functors.
We identify
$\pi\in {\rm Hom}(M,\mathbf{E})(R)$
with a map
$\pi: {\rm Spec}(R)\to {\rm Hom}(M,\mathbf{E})$.

\begin{definition}\rm
For $\pi\in {\rm Hom}(M,\mathbf{E})(R)$,
we define a functor 
\[ {\rm Level}(M,\mathbf{X},\pi)\]
over ${\rm Spec}(R)$
to be the pullback of 
${\rm Level}(M,\mathbf{X})\to{\rm Hom}(M,\mathbf{E})$
along $\pi$.
Since ${\rm Level}(M,\mathbf{X})$
is representable,
so is ${\rm Level}(M,\mathbf{X},\pi)$.
We define $D(M,\mathbf{X},\pi)$
to be the representing ring so that
\[ {\rm Level}(M,\mathbf{X},\pi)=
   {\rm Spec}(D(M,\mathbf{X},\pi)).\]
\end{definition}

By Proposition~\ref{prop:general-level-restriction-height-change},
we obtain the following proposition.

\begin{proposition}
\label{prop:trivilhom-general}
For the zero homomorphism $0\in {\rm Hom}(M,\mathbf{E})(R)$,
we have an isomorphism 
\[ D(M,\mathbf{X},0)\cong 
   D(M,\mathbf{Y}).\]
\end{proposition}

In the following of this subsection
we assume that $\mathbf{E}$ is constant.
Then we have a decomposition
\[ {\rm Level}(M,\mathbf{X})=
   \coprod_{\pi}{\rm Level}(M,\mathbf{X},\pi),\]
where the coproduct ranges over 
$\pi\in {\rm Hom}(M,\mathbf{E}(R))$.
Hence we obtain the following proposition.

\begin{proposition}\label{prop:product-decomposition-DAGp}
If $\mathbf{E}$ is constant,
then there exists a 
decomposition 
\[ D(M,\mathbf{X})\cong  
   \prod_{\pi} D(M,\mathbf{X},\pi),\]
where the product ranges over $\pi\in 
{\rm Hom}(M,\mathbf{E}(R))$.
\end{proposition}

Next, we study the representing ring $D(M,\mathbf{X},\pi)$
for $\pi: M\to \mathbf{E}(R)$.
We set $N=\ker\pi$.

\begin{proposition}\label{prop:pi-component-empty-local}
We assume that $R$ is a local ring with residue field $k$
of characteristic $p$
and that $\mathbf{Y}_k$ is formal.
If {\rm $p$-rank($N$)} is at most the height of $\mathbf{Y}$,
then $D(M,\mathbf{X},\pi)$ is a local ring.
Otherwise,
${\rm Level}(M,\mathbf{X},\pi)=\emptyset$,
and hence $D(M,\mathbf{X},\pi)=0$.
\end{proposition}

\proof
Let $\overline{k}$ be an algebraic closure
of $k$.
It is sufficient to show that
${\rm Level}(M,\mathbf{X},\pi)(\overline{k})$
is an one element set if $p$-rank($N$) is at most 
the height of $\mathbf{Y}$
and that it is empty otherwise. 
Since the map $q$ induces an isomorphism 
$\mathbf{X}(\overline{k})\stackrel{\cong}
{\to}\mathbf{E}(\overline{k})$ of abelian groups,
there is at most one homomorphism
$\phi: M\to \mathbf{X}(\overline{k})$
such that $q\circ \phi=\pi$.
By Proposition~\ref{prop:general-level-restriction-height-change},
$\phi$ is a level $M$-structure on 
$\mathbf{X}_{\overline{k}}$
if and only if 
the restriction $\phi'$ is a level $N$-structure
on $\mathbf{Y}_{\overline{k}}$.
Since $\mathbf{Y}_{\overline{k}}$ is formal,
this is the case if and only if $p$-rank($N$)
is at most 
the height of $\mathbf{Y}_{\overline{k}}$.
\qqq



\subsection{Level structures on $\mathbf{G}_{\mathbb{B}}$}

Recall that $\mathbf{G}=\mathbf{F}_{n+1}[p^{\infty}]$
is the $p$-divisible group associated
to the formal Lie group $\mathbf{F}_{n+1}$ over $E_{n+1}^0$.
We denote by $\mathbf{G}_{\mathbb{B}}$
the base change along the map
${\rm ch}: E_{n+1}^0\to\mathbb{B}_n^0$.
In this subsection we study
level structures on the $p$-divisible group $\mathbf{G}_{\mathbb{B}}$.

First, we recall a connected-\'{e}tale exact sequence
associated to the $p$-divisible group $\mathbf{G}_{\mathbb{B}}$.
Let $\mathbf{H}_{\mathbb{B}}$
be the $p$-divisible group obtained 
from $\mathbf{H}$ by the base change along the map
${\rm inc}: E_n^0\to\mathbb{B}_n^0$,
where $\mathbf{H}$ is the $p$-divisible group
associated to the formal Lie group $\mathbf{F}_n$
over $E_n^0$. 
The maps
${\rm ch}: E_{n+1}^0(\mathbb{CP}^{\infty})\to 
\mathbb{B}_n^0(\mathbb{CP}^{\infty})$ and 
${\rm inc}: E_n^0(\mathbb{CP}^{\infty})\to
\mathbb{B}_n^0(\mathbb{CP}^{\infty})$
induce an isomorphism
\[ \theta: \mathbf{G}_{\mathbb{B}}^0\stackrel{\cong}{\to}
\mathbf{H}_{\mathbb{B}}\]
of connected $p$-divisible groups,
where $\mathbf{G}_{\mathbb{B}}^0$ is the identity component
of $\mathbf{G}_{\mathbb{B}}$.
By \cite[Theorem~5.3]{Torii6},
there is an exact sequence
of $p$-divisible groups
\begin{align}\label{eq:fundamental-exact-squence-Fn-Fn+1-constant}
 0\to \mathbf{H}_{\mathbb{B}}\longrightarrow
        \mathbf{G}_{\mathbb{B}}
        \stackrel{\phantom{q}}{\longrightarrow}
        (\mathbb{Q}_p/\mathbb{Z}_p)_{\mathbb{B}}\to 0
\end{align}
over $\mathbb{B}_n^0$. 

Next, we study the ring
$D(M,\mathbf{G}_{\mathbb{B}})$
for a finite abelian $p$-group $M$.
By Proposition~\ref{prop:level-representability-specR}, 
we have an isomorphism
\[ D(M,\mathbf{G}_{\mathbb{B}})\cong
   {\mathbb{B}_n^0}\otimes_{E_{n+1}^0}D(M,\mathbf{F}_{n+1})\]
of $\mathbb{B}_n^0$-algebras.
By Propositions~\ref{prop:product-decomposition-DAGp} and 
\ref{prop:pi-component-empty-local},
we have a decomposition
\[ D(M,\mathbf{G}_{\mathbb{B}})\cong
   \prod_{\pi} D(M,\mathbf{G}_{\mathbb{B}},\pi),\]
where the product ranges over
$\pi\in {\rm Hom}(M,\mathbb{Q}_p/\mathbb{Z}_p)$
such that $p$-rank$(\ker(\pi))\le n$.

\begin{proposition}
If {\rm $p$-rank$(\ker(\pi))\le n$},
then $D(M,\mathbf{G}_{\mathbb{B}},{\pi})$
is a complete regular local ring.
\end{proposition}

\proof
We set $N=\ker(\pi)$.
For simplicity, we put $D(M)=D(M,\mathbf{G}_{\mathbb{B}},{\pi})$
and $D(N)=D(N,\mathbf{H}_{\mathbb{B}})$.
By Proposition~\ref{prop:pi-component-empty-local},
$D(M)$ is a complete local ring.
If $\pi=0$, then 
$D(M)$ is isomorphic to $D(N)$
by Proposition~\ref{prop:trivilhom-general},
and hence it is regular.
Thus, we assume that $\pi\neq 0$.
Since the Krull dimension of $D(M)$ is $n$,
it is sufficient to show that
the maximal ideal of $D(M)$ is generated by
$n$ elements.

Let $R=E_{n+1}^0/(p,u_1,\ldots,u_{n-1})\cong
\mathbb{F}_{p^{n+1}}\power{u_n}$,
and let $\mathbf{V}$ be the formal group
obtained from $\mathbf{F}_{n+1}$ by
base change along the map $E_{n+1}^0\to R$.
We suppose that $M/N\cong \mathbb{Z}/p^r$.
By the Weierstrass preparation theorem, 
there is a decomposition of the $p^r$-series
$[p^r](x)$ of $\mathbf{V}$ as
$[p^r](x)=g_r(x^{p^{nr}}) u_r(x)$,
where $g_r(y)$ is a monic polynomial in 
$R[y]$ of degree $p^r$, 
and $u_r(x)$ is a unit in $R\power{x}$. 

Let $L$ be the residue field of the local ring $\mathbb{B}_n^0$.
The group scheme $\mathbf{V}[p^r]_L$ 
over $L$
is represented by the ring 
$L[x]\left/\left(g_r\left(x^{p^{nr}}\right)\right)\right.$
and the etale quotient $\mathbf{V}[p^{r}]_L^{\rm et}$
is represented by $L[y]/(g_r(y))$. 
Furthermore, the epimorphism 
$q:\mathbf{V}[p^r]_L\to \mathbf{V}[p^r]^{\rm et}_L$
is represented by a ring homomorphism 
$L[y]/(g_r(y))\to 
L[x]\left/\left(g_r\left(x^{p^{nr}}\right)\right)\right.$ 
given by $y\mapsto x^{p^{nr}}$.

We fix an isomorphism
$\mathbf{V}[p^{\infty}]_L^{\rm et}\cong (\mathbb{Q}_p/\mathbb{Z}_p)_L$,
and take $m\in M$ which is a generator in $M/N$.
Let $U$ be an $L$-algebra, and
let $\phi: M\to \mathbf{V}[p^{\infty}]_L(U)$ be a homomorphism
such that $\phi(N)=0$ and $q\circ\phi=\pi$.
By Proposition~\ref{prop:general-level-restriction-height-change},
we see that $\phi$ is a level $M$-structure on
$\mathbf{V}[p^{\infty}]_U$.
Hence we obtain an isomorphism of functors
\[ {\rm Level}(M,\mathbf{G}_{\mathbb{B}},\pi)
   \times_{{\rm Level}(N,\mathbf{H}_{\mathbb{B}})}
   {\rm Spec}(L)\cong
   \mathbf{V}[p^r]_L\times_{\mathbf{V}[p^r]_L^{\rm et}}
   \{\pi(m)\},\]
where $\{\pi(m)\}$ is the connected component of
$\mathbf{V}[p^{\infty}]_L^{\rm et}\cong (\mathbb{Q}_p/\mathbb{Z}_p)_L$
corresponding to $\pi(m)\in \mathbb{Q}_p/\mathbb{Z}_p$.
This implies an isomorphism of rings
\[ D(M)\otimes_{D(N)} L\cong L[x]
   \left/\left(x^{p^{nr}}-\alpha_r\right)\right. ,\]
where $\alpha_r\in L$ is a root of $g_r(y)$
corresponding to $\{\pi(m)\}$.
By \cite[Lemma~2.7]{Torii6},
we see that $L[x]/\left(x^{p^{n}}-\alpha_r\right)$ is a field.
This implies that the maximal ideal of $D(M)$
is $\mathfrak{m}D(M)$,
where $\mathfrak{m}$ is the maximal ideal of $D(N)$.
Since $D(N)$ is regular of dimension $n$,
$\mathfrak{m}$ is generated by $n$ elements.
Hence the maximal ideal of $D(M)$ is generated by $n$ elements
and $D(M)$ is a regular local ring.
\qqq

\bigskip

Recall that $\Lambda^{n+1}=(\mathbb{Q}_p/\mathbb{Z}_p)^{n+1}$.
For a nonnegative integer $r$,
we set
\[ \begin{array}{rcl}
    {D}_{\mathbb{B}}(r)&=&
   D(\Lambda^{n+1}[p^r],\mathbf{G}_{\mathbb{B}}),\\[2mm]
   {D}_{\mathbb{B}}(r,\pi) &=&
   D(\Lambda^{n+1}[p^r],\mathbf{G}_{\mathbb{B}},\pi),\\
   \end{array}\]
where 
$\pi: \Lambda^{n+1}[p^r]\to \mathbb{Q}_p/\mathbb{Z}_p$
is a homomorphism.
The inclusion $\Lambda^{n+1}[p^r]
\to \Lambda^{n+1}[p^{r+1}]$
induces a $\mathbb{B}_n^0$-algebra homomorphism
$D_{\mathbb{B}}(r)\to D_{\mathbb{B}}(r+1)$, 
and we define a $\mathbb{B}_n^0$-algebra $D_{\mathbb{B}}$ to be
a colimit of $\{D_{\mathbb{B}}(r)\}_{r\ge 0}$:
\[ D_{\mathbb{B}}=\ \subrel{r}{\rm colim}D_{\mathbb{B}}(r).\]
Note that $D_{\mathbb{B}}(r)\cong 
\mathbb{B}_n^0\otimes_{E_{n+1}^0}D_{n+1}(r)$ and
$D_{\mathbb{B}}\cong \mathbb{B}_n^0\otimes_{E_{n+1}^0}D_{n+1}$. 
\if
This implies that the
induced homomorphism $\phi': N\to
\mathbf{F}_n[p^{\infty}](R)$ is not a level $N$-structure,
and that $\phi$ is not a level $M$-structure.
Hence we see that ${\rm Level}_{\pi}
(\Lambda^{n+1}(r),\mathbf{F}_{n+1}[p^{\infty}])=\emptyset$.
\fi

\begin{proposition}
There is an isomorphism
\[ D_{\mathbb{B}}(r)\cong 
   \prod_{\pi} D_{\mathbb{B}}(r,{\pi}),\]
where the product ranges over the 
surjective homomorphisms
$\pi:
\Lambda^{n+1}[p^r]\to (\mathbb{Q}_p/\mathbb{Z}_p)[p^r]$.
\end{proposition}

\proof
The proposition
follows from 
Propositions~\ref{prop:product-decomposition-DAGp} and 
\ref{prop:pi-component-empty-local}
since $p$-rank$(\ker(\pi))>n$
if $\pi$ is not surjective to 
$(\mathbb{Q}_p/\mathbb{Z}_p)[p^r]$.
\qqq


Now, we study an action of
the automorphism group ${\rm Aut}(\Lambda^{n+1}[p^r])$
on $D_{\mathbb{B}}(r)$.
For simplicity, we set
\[ G(r)= {\rm Aut}(\Lambda^{n+1}[p^r]). \]
Note that $G(r)$ is isomorphic to
${\rm GL}_{n+1}(\mathbb{Z}/p^r\mathbb{Z})$ and
acts on the ring ${D}_{\mathbb{B}}(r)$.
Let $P(r,{\pi})$ be the subgroup 
of $G(r)$ 
consisting of the elements
fixing $\pi$:
\[ P(r,{\pi})=\{g\in G(r)|\
            \pi\circ g=\pi\}.\]
We see that the $G(r)$-module
$D_{\mathbb{B}}(r)$ is a coinduced module
of the $P(r,{\pi})$-module $D_{\mathbb{B}}(r,{\pi})$:
\[ D_{\mathbb{B}}(r)\cong\, {\rm Map}_{P(r,{\pi})}({G(r)}, 
   D_{\mathbb{B}}(r,{\pi}))\]
if $\pi:
\Lambda^{n+1}[p^r]\to (\mathbb{Q}_p/\mathbb{Z}_p)[p^r]$
is surjective.

\begin{proposition}\label{prop:invariant-DB}
If $\pi: \Lambda^{n+1}[p^r]\to (\mathbb{Q}_p/\mathbb{Z}_p)[p^r]$ 
is surjective,
then the invariant subring 
of $D_{\mathbb{B}}(r,{\pi})$
under the action of $P(r,\pi)$ is 
$\mathbb{B}_n^0$.
\end{proposition}

\proof
The proposition follows from
\cite[Proposition~6.4]{Torii6}
since we have an isomorphism 
$P(r,\pi)\cong B(W,V)(r)$
and we can identify 
the $P(r,\pi)$-ring $D_{\mathbb{B}}(r,\pi)$ 
with the $B(W,V)(r)$-ring $\mathbb{J}_n(r)$. 
\qqq   

\begin{corollary}\label{cor:invariant-DB}
The invariant ring of $D_{\mathbb{B}}(r)$
under the action of $G(r)$ is 
$\mathbb{B}_n^0$.
\end{corollary}

\proof
The corollary follows from the fact that
the $G(r)$-module
$D_{\mathbb{B}}(r)$ is a coinduced 
module of the $P(r,{\pi})$-module $D_{\mathbb{B}}(r,{\pi})$.
\qqq



\subsection{The action of $\mathbb{G}_{n+1}$ 
on $D(M,\mathbf{G}_{\mathbb{B}})$}

Recall that $\mathbb{G}_{n+1}$ is the $(n+1)$st extended
Morava stabilizer group.
In this subsection 
we study an action of $\mathbb{G}_{n+1}$ 
on the representing ring $D(M,\mathbf{G}_{\mathbb{B}})$
of level structures on the $p$-divisible group
$\mathbf{G}_{\mathbb{B}}$,
where $M$ is a finite abelian $p$-group.

The group $\mathbb{G}_{n+1}$ acts on 
the connected-\'{e}tale exact sequence 
(\ref{eq:fundamental-exact-squence-Fn-Fn+1-constant}),
and hence we obtain a map of exact sequences
\begin{align}\label{eq:map-of-fundamental-exact-seq}
    \begin{array}{ccccccccc}
    0 \to & \mathbf{H}_{\mathbb{B}} & \longrightarrow
    & \mathbf{G}_{\mathbb{B}} & \longrightarrow 
    & (\mathbb{Q}_p/\mathbb{Z}_p)_{\mathbb{B}} & \to & 0\\[1mm]
    & \phantom{\mbox{\scriptsize$\rho_{\mathbf{H}}(g$}}
      \bigg\downarrow\mbox{\scriptsize$\rho_{\mathbf{H}}(g)$} & 
    & \phantom{\mbox{\scriptsize$\rho_{\mathbf{G}}(g)$}}
      \bigg\downarrow\mbox{\scriptsize$\rho_{\mathbf{G}}(g)$} & 
    & \phantom{\mbox{aaaa}}
      \bigg\downarrow
      \mbox{\scriptsize$\rho_{\mathbf{E}}(g)$} & & \\[4mm]
    0 \to & \mathbf{H}_{\mathbb{B}} & \longrightarrow
    & \mathbf{G}_{\mathbb{B}} & \longrightarrow 
    & (\mathbb{Q}_p/\mathbb{Z}_p)_{\mathbb{B}} & \to & 0\\
   \end{array}
\end{align}
for any $g\in \mathbb{G}_{n+1}$ covering
the action on $\mathbb{B}_n^0$.
Note that $\rho_{\mathbf{H}}(g)$ is the identity map 
for any $g\in \mathbb{G}_{n+1}$
(see \cite[\S4]{Torii3}).

Let $M$ be a finite abelian $p$-group.
The action of $\mathbb{G}_{n+1}$ on $E_{n+1}^0$
extends to an action on the complete local ring 
$D(M,\mathbf{F}_{n+1})$.
Since $D(M,\mathbf{G}_{\mathbb{B}})=
\mathbb{B}_n^0\otimes_{E_{n+1}^0}D(M,\mathbf{F}_{n+1})$,
we obtain an action of $\mathbb{G}_{n+1}$ 
on $D(M,\mathbf{G}_{\mathbb{B}})$
by the diagonal action.
We define an action of $\mathbb{G}_{n+1}$
on ${\rm Hom}(M,\mathbb{Q}_p/\mathbb{Z}_p)$ by
\[ g \pi= \rho_{\mathbf{E}}(g)\circ\pi,\] 
where we regard $\rho_{\mathbf{E}}(g)$
as a homomorphism of abelian groups.
We easily obtain the following proposition.

\begin{proposition}\label{prop:permutation-action-Sn+1-DBM}
The action of $\mathbb{G}_{n+1}$ 
on $D(M,\mathbf{G}_{\mathbb{B}})$
induces a local ring homomorphism
\[ g: D(M,\mathbf{G}_{\mathbb{B}},g\pi)
      \longrightarrow
      D(M,\mathbf{G}_{\mathbb{B}},\pi)\] 
for any $g\in\mathbb{G}_{n+1}$.
\end{proposition}

Let $0$ be the zero homomorphism
from $M$ to $\mathbb{Q}_p/\mathbb{Z}_p$.
The action of $\mathbb{G}_{n+1}$ on $D(M,\mathbf{G}_{\mathbb{B}})$
induces an action of $\mathbb{G}_{n+1}$ on 
$D(M,\mathbf{G}_{\mathbb{B}},0)$
by Proposition~\ref{prop:permutation-action-Sn+1-DBM},
and
we have an isomorphism
$D(M,\mathbf{G}_{\mathbb{B}},0)\cong 
   \mathbb{B}_n^0{\otimes}_{E_n^0}D(M,\mathbf{F}_n)$
by Proposition~\ref{prop:trivilhom-general}.
We suppose that $\mathbb{G}_{n+1}$ trivially acts on $D(M,\mathbf{F}_n)$
and diagonally on $\mathbb{B}_n^0{\otimes}_{E_n^0}D(M,\mathbf{F}_n)$.
Then we obtain the following proposition.

\begin{proposition}
\label{prop:trivil-hom-component-DBMII}
The isomorphism
$D(M,\mathbf{G}_{\mathbb{B}},0)\cong 
   \mathbb{B}_n^0{\otimes}_{E_n^0}D(M,\mathbf{F}_n)$
respects the $\mathbb{G}_{n+1}$-actions.
\end{proposition}

\proof
By diagram~(\ref{eq:map-of-fundamental-exact-seq})
and the isomorphism 
$\theta: \mathbf{G}_{\mathbb{B}}^0\stackrel{\cong}{\to}
\mathbf{H}_{\mathbb{B}}$,
we obtain
$\rho_{\mathbf{H}}(g)\circ \theta=\theta\circ \rho_{\mathbf{G}}(g)^0$
for $g\in\mathbb{G}_{n+1}$.
The proposition follows from the fact that
$\rho_{\mathbf{H}}(g)$ is the identity map 
for any $g\in\mathbb{G}_{n+1}$.
\qqq

\subsection{Power operations in $\mathbb{B}_n$-theory}

We would like to construct Hecke operators
in $\mathbb{B}_n$-theory.
For this purpose,
in this subsection
we construct ring operations in $\mathbb{B}_n$-theory
which are extensions of power operations
in $E_{n+1}$-theory. 
We also study compatibility of these operations
with the action of the extended Morava stabilizer
group $\mathbb{G}_{n+1}$.

First, we shall construct a ring operation
\[ 
   \Psi_{M}^{\mathbb{B}}:
   \mathbb{B}_n^0(X)\longrightarrow  
      \mathbb{B}(M)_n^0(X)\]
for a finite abelian $p$-group $M$,
which is an extension of the ring operation
$\Psi_M^{E_{n+1}}: E_{n+1}^0(X)\to E(M)_{n+1}^0(X)$,
where 
$\mathbb{B}(M)_n^*(X)=D(M,\mathbf{G}_{\mathbb{B}})
   \otimes_{\mathbb{B}_n^0}\mathbb{B}_n^*(X)$.

For this purpose,
we consider the composite of 
the operation $\Psi_M^{E_{n+1}}$
with a natural map ${\rm ch}(M): 
E(M)_{n+1}^0(X)\to \mathbb{B}(M)_n^0(X)$,
which is a natural extension of 
${\rm ch}: E_{n+1}^0(X)\to \mathbb{B}_n^0(X)$.
Proposition~\ref{prop:product-decomposition-DAGp}
implies a decomposition
\[ \mathbb{B}(M)_n^*(X)\cong\prod_{\pi}
   \mathbb{B}(M,\pi)^*(X) \]
of generalized cohomology theories,
where 
$\mathbb{B}(M,\pi)_n^*(X)=
   D(M,\mathbf{G}_{\mathbb{B}},\pi)
   {\otimes}_{\mathbb{B}_n^0}\mathbb{B}_n^*(X)$.
We write ${\rm ch}(M,{\pi}):
E(M)_{n+1}^*(X)\to\mathbb{B}(M,\pi)_n^*(X)$
for the composite
$p_{\pi}\circ {\rm ch}(M)$,
where $p_{\pi}: \mathbb{B}(M)_n^*(X)\to
   \mathbb{B}(M,{\pi})_n^*(X)$ 
is the projection.

When $X$ is a one point space,
we obtain a ring homomorphism
${\rm ch}(M,{\pi}): D(M,\mathbf{F}_{n+1})\to 
D(M,\mathbf{G}_{\mathbb{B}},{\pi})$.
The map ${\rm ch}(M,\pi)$ induces
an isogeny of $p$-divisible groups
\[ {\rm ch}(M,\pi)^*\Psi^{E_{n+1}}_M:
   (i^{\mathbb{B}})^*\mathbf{G}_{\mathbb{B}}
   \longrightarrow
   {\rm ch}(M,\pi)^*(\Psi^{E_{n+1}}_M)^*\mathbf{G}\]
from the isogeny of formal groups
$\Psi^{E_{n+1}}_M: i^*\mathbf{F}_{n+1}\to
(\Psi^{E_{n+1}}_M)^*\mathbf{F}_{n+1}$,
where $i^{\mathbb{B}}: \mathbb{B}_n^0\to
D(M,\mathbf{G}_{\mathbb{B}},\pi)$ is the inclusion map.

Let $\phi: M\to \mathbf{G}_{\mathbb{B}}(R)$
be a level $M$-structure,
which corresponds to a map 
$D(M,\mathbf{G}_{\mathbb{B}},\pi)\to R$ in $\mathcal{CL}$.
By exact 
sequence~(\ref{eq:fundamental-exact-squence-Fn-Fn+1-constant}),
we obtain a homomorphism $\phi': N\to \mathbf{H}_{\mathbb{B}}(R)$,
where $N=\ker(\pi)$. 
By Proposition~\ref{prop:general-level-restriction-height-change},
$\phi'$ is a level $N$-structure on $\mathbf{H}_R$,
and hence we obtain a local ring homomorphism
${\rm inc}(M,\pi): D(N,\mathbf{F}_n)\to
              D(M,\mathbf{G}_{\mathbb{B}},\pi)$,
which induces a multiplicative stable
cohomology operation
\[ {\rm inc}(M,\pi): E(N)_n^*(X)\longrightarrow
                    \mathbb{B}(M,\pi)_n^*(X).\]  

The following lemma 
follows from Proposition~\ref{prop:exact-seq-quotient-pdivisible}.

\begin{lemma}\label{lemma:iso-identity-FM-FN}
The isomorphism
$\theta: \mathbf{G}_{\mathbb{B}}^0\stackrel{\cong}{\to}
\mathbf{H}_{\mathbb{B}}$ induces
an isomorphism 
\[ \overline{\theta}: (\mathbf{G}_R/[\phi(M)])^0
   \stackrel{\cong}{\longrightarrow}
   \mathbf{H}_R/[\phi'(N)]\]
of connected $p$-divisible groups over $R$,
where $(\mathbf{G}_R/[\phi(M)])^0$ is the identity
component of the quotient $p$-divisible group
$\mathbf{G}_R/[\phi(M)]$.
\end{lemma}


Using Lemma~\ref{lemma:iso-identity-FM-FN},
we shall show that
${\rm ch}(M,\pi)\circ \Psi_M^{E_{n+1}}:
E_{n+1}^0\to D(M,\mathbf{G}_{\mathbb{B}},\pi)$
extends to a map 
$f: \mathbb{A}_n^0\to D(M,\mathbf{G}_{\mathbb{B}},\pi)$
of local rings.

\begin{lemma}\label{lemma:extension-chMpi-psiMen+1}
The composite
${\rm ch}(M,\pi)\circ \Psi_M^{E_{n+1}}:
E_{n+1}^0\to D(M,\mathbf{G}_{\mathbb{B}},\pi)$
extends to a map 
$f: \mathbb{A}_n^0\to D(M,\mathbf{G}_{\mathbb{B}},\pi)$
of local rings.
\end{lemma}

\proof
Put $R=D(M,\mathbf{G}_{\mathbb{B}},\pi)$
and $\psi={\rm ch}(M,\pi)\circ \Psi_M^{E_{n+1}}$.
We have an isomorphism
$\psi^*\mathbf{G}\cong\mathbf{G}_R/[\phi(M)]$
of $p$-divisible groups.
Lemma~\ref{lemma:iso-identity-FM-FN}
implies an isomorphism
$(\psi^*\mathbf{G})^0\cong \mathbf{H}_R/[\phi'(N)]$.
Since the height of $\mathbf{H}_R$ is $n$,
we see that the height of $(\psi^*\mathbf{G})^0$ 
is also $n$.
This implies 
$\psi(u_i)=0$ for $i<n$, and $\psi(u_n)\neq 0$
in the residue field of $R$.
Hence $\psi$ extends to a map 
$f: \mathbb{A}_n^0\to R$ 
since $R$ is a complete local ring.
\qqq

\begin{proposition}\label{prop:relations-Psi-ch-i}
There exists a unique map of local rings
\[ 
   \Psi_{M,\pi}^{\mathbb{B}}:
   \mathbb{B}_n^0\longrightarrow 
   D(M,\mathbf{G}_{\mathbb{B}},\pi) \]
satisfying 
\[ \begin{array}{lcl}
    \Psi_{M,\pi}^{\mathbb{B}}\circ {\rm ch}&=&
    {\rm ch}(M,{\pi})\circ\Psi_{M}^{E_{n+1}},\\[2mm]
    \Psi_{M,\pi}^{\mathbb{B}}\circ {\rm inc}&=&
    {\rm inc}(M,\pi)\circ \Psi_N^{E_n},
   \end{array}\]
and making the following diagram commute
\[ \begin{array}{ccc}
    (i^{\mathbb{B}})^*\mathbf{G}_{\mathbb{B}}^0 &
    \stackrel{{\rm ch}(M,\pi)^*\Psi^{E_{n+1}}_M}
    {\hbox to 20mm{\rightarrowfill}}&
    (\Psi^{\mathbb{B}}_{M,\pi})^*\mathbf{G}_{\mathbb{B}}^0\\[2mm]
    \mbox{$\scriptstyle (i^{\mathbb{B}})^*\theta$}
    \bigg\downarrow
    \phantom{\mbox{$\scriptstyle (i^{\mathbb{B}})^*\theta$}}
    &&
    \phantom{\mbox{$\scriptstyle (\Psi^{\mathbb{B}}_{M,\pi})^*\theta$}}
    \bigg\downarrow
    \mbox{$\scriptstyle (\Psi^{\mathbb{B}}_{M,\pi})^*\theta$}\\[2mm]
    (i^{\mathbb{B}})\mathbf{H}_{\mathbb{B}}&
    \stackrel{{\rm inc}(M,\pi)^*\Psi^{E_n}_N}
    {\hbox to 20mm{\rightarrowfill}}&
    (\Psi^{\mathbb{B}}_{M,\pi})^*\mathbf{H}_{\mathbb{B}}.
   \end{array}\]
\end{proposition}

\proof
Put $R=D(M,\mathbf{G}_{\mathbb{B}},\pi)$.
The isogeny
${\rm ch}(M,\pi)^*\Psi^{E_{n+1}}_M:
(i^{\mathbb{B}})^*\mathbf{G}_{\mathbb{B}}=
\mathbf{G}_R
\to f^*\mathbf{G}_{\mathbb{A}}$
of $p$-divisible groups 
induces an isomorphism
$f^*\mathbf{G}_{\mathbb{A}}\cong \mathbf{G}_R/[\phi(M)]$,
where $f: \mathbb{A}_n^0\to R$ 
is the extension of ${\rm ch}(M,\pi)\circ\Psi_M^{E_{n+1}}$
given in Lemma~\ref{lemma:extension-chMpi-psiMen+1}.
We also have the isogeny
${\rm inc}(M,\pi)^*\Psi^{E_n}_N:
(i^{\mathbb{B}})\mathbf{H}_{\mathbb{B}}=
\mathbf{H}_R\to g^*\mathbf{H}$,
which induces an isomorphism
$g^*\mathbf{H}\cong \mathbf{H}_R/[\phi'(N)]$,
where $g={\rm inc}(M,\pi)\circ \Psi_N^{E_n}$. 
By Lemma~\ref{lemma:MUPMUP-variation},
it suffices to show that 
there is an isomorphism
$\overline{\theta}: f^*\mathbf{G}_{\mathbb{A}}^0
\cong
g^*\mathbf{H}$  
of $p$-divisible groups
over $D(M,\mathbf{G}_{\mathbb{B}},\pi)$.
We obtain the desired isomorphism
by Lemma~\ref{lemma:iso-identity-FM-FN}. 
\qqq


\bigskip

When $X$ is a finite complex,
we have a natural isomorphism
$\mathbb{B}_n^0(X)\cong \mathbb{B}_n^0
   \otimes_{E_{n+1}^*}E_{n+1}^0(X)$.
We define a natural map
\[   
   \Psi_{M,\pi}^{\mathbb{B}}:
   \mathbb{B}_n^0(X)\longrightarrow \mathbb{B}(M,\pi)_n^0(X) \]
to be the composition
\[ \begin{array}{rcl}
   \mathbb{B}_n^0\subrel{E_{n+1}^0}{\otimes}E_{n+1}^0(X)&
   \stackrel{\Psi_{M,\pi}^{\mathbb{B}}\otimes
   \psi}{\hbox to 14mm{\rightarrowfill}}&
   D(M,\mathbf{G}_{\mathbb{B}},\pi)
   \subrel{\psi,E_{n+1}^0,\psi}{\otimes}
   D(M,\mathbf{G}_{\mathbb{B}},\pi)
   \subrel{E_{n+1}^0}{\otimes} 
   E_{n+1}^0(X)\\[2mm]
   &\stackrel{m\otimes 1}{\hbox to 14mm{\rightarrowfill}}&
   D(M,\mathbf{G}_{\mathbb{B}},{\pi})
   \subrel{E_{n+1}^0}{\otimes}
   E_{n+1}^0(X), 
   \end{array}\]
where $\psi={\rm ch}(M,\pi)\circ \Psi_M^{E_{n+1}}$ and
$m$ is the multiplication map of 
$D(M,\mathbf{G}_{\mathbb{B}},{\pi})$.

When $X$ is a CW-complex,
we let $\{X_{\alpha}\}$ be the filtered system 
of all finite subcomplexes of $X$.  
Since $\mathbb{B}_n^0$ is a complete Noetherian local ring,
we have an isomorphism
$\mathbb{B}_n^0(X)\cong\ \subrel{\alpha}{\lim}
                              \mathbb{B}_n^0(X_{\alpha})$.
Since $D(M,\mathbf{G}_{\mathbb{B}},\pi)$ 
is a finitely generated free
$\mathbb{B}_n^0$-module,
we also have an isomorphism
$\mathbb{B}(M,\pi)_n^0(X)\cong\ \subrel{\alpha}{\lim}
    \mathbb{B}(M,\pi)_n^0(X_{\alpha})$.
Hence we see that 
the natural map 
$\Psi_{M,\pi}^{\mathbb{B}}$
uniquely extends to a ring operation 
\[ 
   \Psi_{M,\pi}^{\mathbb{B}}:
   \mathbb{B}_n^0(X)
   \longrightarrow 
   \mathbb{B}(M,\pi)_n^0(X)\]
for any space $X$.

\begin{remark}\label{remark:Psi-M-pi-ch}
\rm
By construction,
we have
$\Psi_{M,\pi}^{\mathbb{B}}\circ {\rm ch}=
{\rm ch}(M,\pi)\circ \Psi_M^{E_{n+1}}$.
\end{remark}

The following theorem
describes a relationship between
the operations $\Psi_{M,\pi}^{\mathbb{B}}$
and $\Psi_N^{E_n}$,
where $N$ is the kernel of
$\pi: M\to \mathbb{Q}_p/\mathbb{Z}_p$. 

\begin{theorem}\label{thm:relation-psi-M-pi-Psi-N}
There is a natural commutative diagram
\[ \begin{array}{ccc}
    E_n^0(X) & 
    \stackrel{\Psi_N^{E_n}}{\hbox to 10mm{\rightarrowfill}}&
    E(N)_n^0(X)\\[1mm]
    \mbox{\scriptsize${\rm inc}$}
    \bigg\downarrow\phantom{\mbox{\scriptsize${\rm inc}$}} 
    & &
    \phantom{\mbox{\scriptsize${\rm inc}(M,\pi)$}}
    \bigg\downarrow
    \mbox{\scriptsize${\rm inc}(M,\pi)$}\\[1mm]
    \mathbb{B}_n^0(X) &
    \stackrel{\Psi_{M,\pi}^{\mathbb{B}}}
    {\hbox to 10mm{\rightarrowfill}} &
    \mathbb{B}_n(M,\pi)_n^0(X)\\ 
   \end{array}\]
for any space $X$.
\end{theorem}

\proof
By \cite[Proposition~3.7]{Ando2},
an unstable ring operation between Landweber exact cohomology theories
is determined by its values on a one point space
and the infinite complex projective space 
$\mathbb{CP}^{\infty}$.
Hence it is sufficient to show that
the diagram is commutative when $X$ is a one point space
and when $X$ is $\mathbb{CP}^{\infty}$.
If $X$ is a one point space,
then the diagram is commutative 
by Proposition~\ref{prop:relations-Psi-ch-i}.
If $X=\mathbb{CP}^{\infty}$, 
then the diagram is commutative 
by Lemma~\ref{lemma:iso-identity-FM-FN}.
\qqq

Assembling the operations
$\Psi_{M,\pi}^{\mathbb{B}}$
for all homomorphisms $\pi: M\to \mathbb{Q}_p/\mathbb{Z}_p$,
we shall define a ring operation 
$\Psi_M^{\mathbb{B}}$ with values
in $\mathbb{B}(M)_n^0(X)$.

\begin{definition}\rm
We define a ring operation 
\[ 
   \Psi_{M}^{\mathbb{B}}:
   \mathbb{B}_n^0(X)
   \longrightarrow 
   \mathbb{B}(M)_n^0(X)\]
to be the one satisfying
$\Psi_{M,\pi}^{\mathbb{B}}=
   (p_{\pi}\otimes 1)\circ \Psi_M^{\mathbb{B}}$
for any homomorphism $\pi: M\to \mathbb{Q}_p/\mathbb{Z}_p$.
\end{definition}

The following proposition is easily obtained
by the definition of $\Psi_M^{\mathbb{B}}$
and Remark~\ref{remark:Psi-M-pi-ch},
which describes a relationship between
the operations $\Psi_M^{\mathbb{B}}$
and $\Psi_M^{E_{n+1}}$.

\begin{proposition}\label{prop-on-Psi-operation-B-M}
We have 
$\Psi_M^{\mathbb{B}}\circ {\rm ch}
   = {\rm ch}(M)\circ \Psi_M^{E_{n+1}}$.
\end{proposition}

Next, we study 
compatibility of the operation
$\Psi_M^{\mathbb{B}}$
with an action of the extended
Morava stabilizer group $\mathbb{G}_{n+1}$.
We suppose that $\mathbb{G}_{n+1}$ diagonally
acts on $\mathbb{B}(M)_n^0(X)= 
D(M,\mathbf{G}_{\mathbb{B}})\otimes_{\mathbb{B}_n^0}
\mathbb{B}_n^0(X)$ when $X$ is a finite complex.
This action uniquely extends to
an action on $\mathbb{B}(M)_n^0(X)$
for any space $X$
since $\mathbb{B}(M)_n^0(X)\cong\,
\subrel{\alpha}{\rm lim}\,\mathbb{B}(M)_n^0(X_{\alpha})$.

\begin{theorem}\label{thm:equivariance-operation-psi-B}
The operation
$\Psi_{M}^{\mathbb{B}}: \mathbb{B}_n^0(X)\to
\mathbb{B}(M)_n^0(X)$
is $\mathbb{G}_{n+1}$-equivariant.
\end{theorem}

In order to prove Theorem~\ref{thm:equivariance-operation-psi-B},
we first study the action of $\mathbb{G}_{n+1}$
on the ring homomorphism
$\Psi_M^{\mathbb{B}}:\mathbb{B}_n^0\to D(M,\mathbf{G}_{\mathbb{B}})
\cong \prod_{\pi}D(M,\mathbf{G}_{\mathbb{B}},\pi)$.
By Lemma~\ref{lemma:MUPMUP-variation},
the local ring homomorphism $\Psi_{M,\pi}^{\mathbb{B}}$ 
is characterized by 
compositions $F=\Psi_{M,\pi}^{\mathbb{B}}\circ{\rm j}$ and
$G=\Psi_{M,\pi}^{\mathbb{B}}\circ {\rm inc}$,
and an isomorphism of formal groups between 
$F^*\mathbf{G}_{\mathbb{A}}^0$
and $G^*\mathbf{H}$.

By Proposition~\ref{prop:relations-Psi-ch-i},
we have $\Psi_{M,\pi}^{\mathbb{B}}\circ{\rm inc}=
{\rm inc}(M,\pi)\circ\Psi_{N}^{E_{n+1}}$,
where $N={\rm ker}(\pi)$.
We consider the action of $\mathbb{G}_{n+1}$
on the ring homomorphisms
${\rm inc}(M,\pi): D(N,\mathbf{F}_n)\to
D(M,\mathbf{G}_{\mathbb{B}},\pi)$.
Here we suppose that $\mathbb{G}_{n+1}$ 
trivially acts on $D(N,\mathbf{F}_n)$ .
The following lemma follows from the fact that
$\rho_{\mathbf{H}}(g)$ in (\ref{eq:map-of-fundamental-exact-seq})
is the identity map
for any $g\in \mathbb{G}_{n+1}$.

\begin{lemma}\label{lemma:commutativity-DN-GN+1-RMN-DBPI}
We have $g\circ {\rm inc}(M,g\pi)={\rm inc}(M,\pi)
: D(N,\mathbf{F}_n)\to D(M,\mathbf{G}_{\mathbb{B}},\pi)$
for any $g\in\mathbb{G}_{n+1}$.
\end{lemma}

Using 
Proposition~\ref{prop:g-Psi-commutation-En},
Lemma~\ref{lemma:MUPMUP-variation},
Proposition~\ref{prop:permutation-action-Sn+1-DBM},
Proposition~\ref{prop:relations-Psi-ch-i},
and
Lemma~\ref{lemma:commutativity-DN-GN+1-RMN-DBPI},
we obtain the following proposition.

\begin{proposition}\label{prop:Bn-PSI-MPI-G-commutativity}
There is a commutative diagram
\[ \begin{array}{ccc}
    \mathbb{B}_n^0 & \stackrel{g}{\hbox to 10mm{\rightarrowfill}} &
    \mathbb{B}_n^0 \\[1mm]
    \mbox{\scriptsize
    $\Psi_{M,g\pi}^{\mathbb{B}}$}
    \bigg\downarrow
    \phantom{\mbox{\scriptsize
    $\Psi_{M,g\pi}^{\mathbb{M}}$}} 
    & & 
    \phantom{\mbox{\scriptsize$\Psi_{M,\pi}^{\mathbb{B}}$}}
    \bigg\downarrow
    \mbox{\scriptsize$\Psi_{M,\pi}^{\mathbb{B}}$} \\[4mm]
    D(M,\mathbf{G}_{\mathbb{B}},{g\pi}) &
    \stackrel{g}{\hbox to 10mm{\rightarrowfill}} &
    D(M,\mathbf{G}_{\mathbb{B}},{\pi})\\[1mm]
   \end{array}\]
for any $g\in \mathbb{G}_{n+1}$.
\end{proposition}

\proof
Put $R=D(M,\mathbf{G}_{\mathbb{B}},\pi)$.
By 
Propositions~\ref{prop:g-Psi-commutation-En},
\ref{prop:permutation-action-Sn+1-DBM},
and \ref{prop:relations-Psi-ch-i},
we have
$g\circ \Psi_{M,g\pi}^{\mathbb{B}}\circ {\rm ch}
=\Psi_{M,\pi}^{\mathbb{B}}\circ g\circ {\rm ch}$.
This induces a map $F_g: \mathbb{A}_n^0\to R$ in $\mathcal{CL}$,
and there are canonical isomorphisms
$F_g^*\mathbf{G}_{\mathbb{A}}^0\cong
(g^*\mathbf{G}_R/g^*M)^0\cong g^*(\mathbf{G}_R/M)^0$.
On the other hand,
by 
Proposition~\ref{prop:relations-Psi-ch-i} and
Lemma~\ref{lemma:commutativity-DN-GN+1-RMN-DBPI},
$g\circ \Psi_{M,g\pi}^{\mathbb{B}}\circ {\rm inc}
=\Psi_{M,\pi}^{\mathbb{B}}\circ g\circ {\rm inc}$.
We denote this map by $G_g$.
There are canonical isomorphisms
$G_g^*\mathbf{H}\cong
g^*\mathbf{H}_R/g^*N\cong g^*(\mathbf{H}_R/N)$,
where $N=\ker(\pi)$.
The isomorphism $\theta:\mathbf{G}_{\mathbb{B}}^0
\stackrel{\cong}{\to}\mathbf{H}_{\mathbb{B}}$
induces isomorphisms
$\overline{g^*\theta}: (g^*\mathbf{G}_R/g^*M)^0\stackrel{\cong}{\to}
g^*\mathbf{H}_R/g^*N$ and
$g^*\overline{\theta}: g^*(\mathbf{G}_R/M)^0\stackrel{\cong}{\to}
g^*(\mathbf{H}_R/N)$.
The proposition follows from
Lemma~\ref{lemma:MUPMUP-variation}
since we have $\overline{g^*\theta}=g^*\overline{\theta}$
under the above isomorphisms.
\qqq

\begin{corollary}\label{cor:equivariance-psi-B-one-point-case}
The map $\Psi_M^{\mathbb{B}}: \mathbb{B}_n^0\to 
D(M,\mathbf{G}_{\mathbb{B}})$ 
is $\mathbb{G}_{n+1}$-equivariant.
\end{corollary}

\proof[Proof of Theorem~\ref{thm:equivariance-operation-psi-B}]
It is sufficient to show that 
the theorem holds for any finite complex $X$.
In this case,
since $\Psi_M^{\mathbb{B}}$ is a ring operation and
$\mathbb{B}_n^0(X)\cong 
\mathbb{B}_n^0\otimes_{E_{n+1}^0}E_{n+1}^0(X)$,
it is sufficient to show that this holds for
a one point space and that $\Psi_M^{\mathbb{B}}\circ{\rm ch}$
is $\mathbb{G}_{n+1}$-equivariant.
The case when $X$ is a one point space follows
from Corollary~\ref{cor:equivariance-psi-B-one-point-case}.
The equivariance of $\Psi_M^{\mathbb{B}}\circ{\rm ch}$
follows from 
Proposition~\ref{prop-on-Psi-operation-B-M},
Proposition~\ref{prop:g-Psi-commutation-En},
and the fact that
${\rm ch}(M)$ is $\mathbb{G}_{n+1}$-equivariant.
\qqq



Now,
we introduce
a ring operation $\Phi_M^{\mathbb{B}}$
which is an extension of
the ring operation
$\Phi_M^{E_{n+1}}$.
We suppose that
$M$ is a subgroup of $\Lambda^{n+1}$.
We define an operation
\[ 
   \Phi_{M}^{\mathbb{B}}:
   \mathbb{B}_n^0(X)\longrightarrow
   D_{\mathbb{B}}\otimes_{\mathbb{B}_n^0}
   \mathbb{B}_n^0(X)\]
by the composition
$\Phi_{M}^{\mathbb{B}}=(I_M^{\mathbb{B}}\otimes 1)\circ 
   \Psi_M^{\mathbb{B}}$,
where
$I_M^{\mathbb{B}}$ 
is the map 
$D(M,\mathbf{G}_{\mathbb{B}})
   \to
   D_{\mathbb{B}}$
induced by the inclusion
$M\hookrightarrow\Lambda^{n+1}$.

Since $D_{\mathbb{B}}\cong \mathbb{B}_n^0\otimes_{E_{n+1}^0}D_{n+1}$,
the map ${\rm ch}: E_{n+1}^0(X)\to \mathbb{B}_n^0(X)$
extends to a map
\[ {\rm ch}(\Lambda^{n+1}): 
    D_{n+1}{\otimes}_{E_{n+1}^0}E_{n+1}^0(X)
    \longrightarrow 
    D_{\mathbb{B}}{\otimes}_{\mathbb{B}_n^0}
    \mathbb{B}_n^0(X).\]

By Lemma~\ref{lemma:restriction-g-level-str-commutativity},
Proposition~\ref{prop-on-Psi-operation-B-M},
and Theorem~\ref{thm:equivariance-operation-psi-B},
we obtain the following proposition.

\begin{proposition}\label{prop:on-Phi-M-ch-and-equivariance}
The operation
$\Phi_{M}^{\mathbb{B}}: \mathbb{B}_n^0(X)\to 
D_{\mathbb{B}}\otimes_{\mathbb{B}}\mathbb{B}_n^0(X)$
is $\mathbb{G}_{n+1}$-equivariant,
and we have
$\Phi_M^{\mathbb{B}}\circ{\rm ch}=
   {\rm ch}(\Lambda^{n+1})\circ\Phi_M^{E_{n+1}}$.
\end{proposition}

In order to compare Hecke operators
in $E_{n+1}$-theory with those in $E_n$-theory
in \S\ref{section:comparison-Hecke-operators},
we introduce a ring $D_{\mathbb{B}}(\mu)$
for a split surjection
$\mu: \Lambda^{n+1}\to \mathbb{Q}_p/\mathbb{Z}_p$,
and a ring operation
$\Phi_{M,\mu}^{\mathbb{B}}:
\mathbb{B}_n^0(X)\to D_{\mathbb{B}}(\mu)\otimes_{\mathbb{B}_n^0}
\mathbb{B}_n^0(X)$
for a subgroup $M$ of $\Lambda^{n+1}$.
We describe a relationship between
the operations $\Phi_{M,\mu}^{\mathbb{B}}$
and $\Phi_N^{E_n}$,
where $N=M\cap \ker(\mu)$.

The inclusion $\Lambda^{n+1}[p^r]\to\Lambda^{n+1}[p^{r+1}]$
of abelian groups induces a map
$D_{\mathbb{B}}(r,\mu[p^r])\to D_{\mathbb{B}}(r+1,\mu[p^{r+1}])$
of local rings.
We set
\[ D_{\mathbb{B}}(\mu)=\ \subrel{r}{\rm colim}
   D_{\mathbb{B}}(r,\mu[p^r]) .\]
Assembling
the projections $D_\mathbb{B}(r)\to D_{\mathbb{B}}(r,\mu[p^r])$
for $r\ge 0$,
we obtain a map 
$p_{\mu}: D_{\mathbb{B}}\to D_{\mathbb{B}}(\mu)$.
We define a ring operation
\[ 
   \Phi_{M,\mu}^{\mathbb{B}}:
   \mathbb{B}_n^0(X)\longrightarrow
   D_{\mathbb{B}}(\mu){\otimes}_{\mathbb{B}_n^0}
   \mathbb{B}_n^0(X)\]
by the composition 
$\Phi_{M,\mu}^{\mathbb{B}}=(p_{\mu}\otimes 1)
   \circ\Phi_{M}^{\mathbb{B}}$.

To describe
the composite $\Phi_{M,\mu}^{\mathbb{B}}\circ{\rm inc}:
E_n^0(X)\to \mathbb{B}_n^0(X)\to
D_{\mathbb{B}}(\mu)\otimes_{\mathbb{B}_n^0}
\mathbb{B}_n^0(X)$,
we set
\[ D_n(V)=\ \subrel{r}{\rm colim} D(V[p^r],\mathbf{F}_n).\]
where $V=\ker(\mu)$.
Since $V\cong \Lambda^n$,
we have an isomorphism $D_n(V)\cong D_n$,
and an operation
$\Phi_N^{E_n}: E_n^0(X)\to D_n(V)\otimes_{E_n^0}E_n^0(X)$,
where $N=M\cap V$.

By Theorem~\ref{thm:relation-psi-M-pi-Psi-N},
we obtain the following proposition
which describes a relationship between
the operations $\Phi_{M,\mu}^{\mathbb{B}}$
and $\Phi_N^{E_n}$.

\begin{proposition}\label{prop:relation-PhiB-PhiEn}
There is a natural commutative diagram
\[ \begin{array}{ccc}
     E_n^0(X) & \stackrel{\Phi_N^{E_n}}
     {\hbox to 10mm{\rightarrowfill}} &
     D_n(V){\otimes}_{E_n^0} E_n^0(X)\\[4mm]
     \mbox{\scriptsize${\rm inc}$}\bigg\downarrow
     \phantom{\mbox{\scriptsize${\rm i}$}} &&
     \phantom{\mbox{\scriptsize${\rm inc}(\mu)
        \otimes {\rm inc}$}}\bigg\downarrow
     \mbox{\scriptsize${\rm inc}(\mu)\otimes {\rm inc}$}\\[2mm]
    \mathbb{B}_n^0(X) &
    \stackrel{\Phi_{M,\mu}^{\mathbb{B}}}
    {\hbox to 10mm{\rightarrowfill}} &
    D_{\mathbb{B}}(\mu){\otimes}_{\mathbb{B}_n^0}
    \mathbb{B}_n^0(X)\\[2mm]
   \end{array}\]
for any space $X$,
where 
${\rm inc}(\mu): D_n(V)\to D_{\mathbb{B}}(\mu)$
is a map induced by
the inclusion $V\hookrightarrow\Lambda^{n+1}$.
\end{proposition}

\subsection{Hecke operators in $\mathbb{B}_n$-theory}

In this subsection
we define Hecke operators in $\mathbb{B}_n$-theory.

Let $M$ be a finite abelian $p$-group.
Recall that $m_{n+1}(M)$ 
is the set of all subgroups of $\Lambda^{n+1}$
which are isomorphic to $M$. 
We consider a natural map given by
\[ \sum_{M'\in m_{n+1}(M)} 
   \Phi_{M'}^{\mathbb{B}}: \mathbb{B}_n^0(X)\longrightarrow
             D_{\mathbb{B}}{\otimes}_{\mathbb{B}_n^0}
                 \mathbb{B}_n^0(X).\]
The image of this map is invariant 
under the action of ${\rm GL}_{n+1}(\mathbb{Z}_p)$.
Since the invariant subring of $D_{\mathbb{B}}$
under the action of ${\rm GL}_{n+1}(\mathbb{Z}_p)$
is $\mathbb{B}_n^0$
by Corollary~\ref{cor:invariant-DB},
we obtain the following theorem
in the same way as in 
the Hecke operators in Morava $E$-theory.

\begin{theorem}\label{thm:Btheory-Hecke-op}
For an abelian $p$-group $M$,
there is an unstable additive cohomology operation
\[ 
   \widetilde{\rm T}_M^{\mathbb{B}}: 
   \mathbb{B}_n^0(X)
   \longrightarrow \mathbb{B}_n^0(X) \]
given by
$(i^{\mathbb{B}}\otimes 1)\circ \widetilde{\rm T}_M^{\mathbb{B}}=
   \sum_{M'\in m_{n+1}(M)} 
   \Phi_{M'}^{\mathbb{B}}$,
where $i^{\mathbb{B}}$ is the $\mathbb{B}_n^0$-algebra structure map
of $D_{\mathbb{B}}$.
For a nonnegative integer $k$,
there is also an unstable additive operation 
\[ 
   \widetilde{\rm T}_{\mathbb{B}}\left(p^k\right):
   \mathbb{B}_n^0(X)
   \longrightarrow \mathbb{B}_n^0(X)\]
given by
$(i^{\mathbb{B}}\otimes 1)\circ
   \widetilde{\rm T}_{\mathbb{B}}\left(p^k\right)=
   \sum_{|M|=p^k}    \Phi_M^{\mathbb{B}}$,
where the sum ranges over all subgroups $M$ of $\Lambda^{n+1}$
such that the order $|M|$ is $p^k$.
\end{theorem}

By
Proposition~\ref{prop:on-Phi-M-ch-and-equivariance},
we obtain the following proposition.

\begin{proposition}\label{prop:Hecke-G_{n+1}-equivariance}
The operations 
$\widetilde{\rm T}_M^{\mathbb{B}}$ and 
$\widetilde{\rm T}_{\mathbb{B}}\left(p^k\right)$
are $\mathbb{G}_{n+1}$-equivariant.
\end{proposition}


Now, we define a natural action of the Hecke algebra
$\mathcal{H}_{n+1}$ on $\mathbb{B}_n^0(X)$
for any space $X$.
Recall that $\Delta_{n+1}={\rm End}(\Lambda^{n+1})\cap
{\rm GL}_{n+1}(\mathbb{Q}_p)$ and
$\Gamma_{n+1}={\rm GL}_{n+1}(\mathbb{Z}_p)$.
The Hecke algebra $\mathcal{H}_{n+1}$
is the endomorphism ring 
of the $\mathbb{Z}[\Delta_{n+1}]$-module
$\mathbb{Z}[\Delta_{n+1}/\Gamma_{n+1}]$. 
There is an additive isomorphism
$\mathcal{H}_{n+1}\cong
\mathbb{Z}[\Gamma_{n+1}\backslash\Delta_{n+1}/\Gamma_{n+1}]$,
and  
we identify $\Gamma_{n+1}\backslash\Delta_{n+1}/\Gamma_{n+1}$
with the set of all isomorphism classes of
finite abelian $p$-groups with $p$-rank $\le n+1$.
For a finite abelian $p$-group $M$ with $p$-rank $\le n+1$,
we denote by $\widetilde{M}$ 
the associated endomorphism
of $\mathbb{Z}[\Delta_{n+1}]$-module
$\mathbb{Z}[\Delta_{n+1}/\Gamma_{n+1}]$.
 
\begin{theorem}
\label{thm:Hn+1-module-structure-on-BnX}
Assigning to $\widetilde{M}$
the Hecke operator $\widetilde{\rm T}_M^{\mathbb{B}}$,
there is a natural $\mathcal{H}_{n+1}$-module
structure on $\mathbb{B}_n^0(X)$
for any space $X$
such that ${\rm ch}: E_{n+1}^0(X)\to \mathbb{B}_n^0(X)$
is a map of $\mathcal{H}_{n+1}$-modules.
\end{theorem}

\proof
By \cite[Lemma~14.4]{Rezk},
there is a ring homomorphism
$\psi_x: D_{n+1}\to D_{n+1}$
for each $x\in\Delta_{n+1}$
such that $\psi_x\circ i=\Phi_M^{E_{n+1}}$,
where $M$ is the kernel of $x:\Lambda^{n+1}\to\Lambda^{n+1}$.
Using the isomorphism
$D_{\mathbb{B}}\otimes_{\mathbb{B}_n^0}\mathbb{B}_n^0(X)\cong
D_{n+1}\otimes_{E_{n+1}^0}\mathbb{B}_n^0(X)$
and
Proposition~\ref{prop:on-Phi-M-ch-and-equivariance},
we can extend $\Phi_M^{\mathbb{B}}$ to a ring operation
\[ \psi_x^{\mathbb{B}}:
   D_{\mathbb{B}}\otimes_{\mathbb{B}_n^0}\mathbb{B}_n^0(X)
   \longrightarrow
   D_{\mathbb{B}}\otimes_{\mathbb{B}_n^0}\mathbb{B}_n^0(X)\]
such that 
$\psi_x^{\mathbb{B}}\circ (i^{\mathbb{B}}\otimes 1)=\Phi_M^{\mathbb{B}}$.
We obtain the theorem
in the same way as in \cite[Proposition~14.3]{Rezk}.
\qqq



By Theorem~\ref{thm:Hn+1-module-structure-on-BnX},
there is a natural $\mathcal{H}_{n+1}$-module 
structure on $\mathbb{B}_n^0(X)$ for any space $X$.
By Proposition~\ref{prop:Hecke-G_{n+1}-equivariance},
the $\mathcal{H}_{n+1}$-module structure
commutes with the action of $\mathbb{G}_{n+1}$.
This implies an $\mathcal{H}_{n+1}$-module
structure on the invariant submodule.
In \cite{Torii3} we showed that
the map ${\rm inc}$ induces
a natural isomorphism
\[ E_n^0(X)\stackrel{\cong}{\to}
   (\mathbb{B}_n^0(X))^{\mathbb{G}_{n+1}} \]
for any spectrum $X$,
where the right hand side is the $\mathbb{G}_{n+1}$-invariant
submodule of $\mathbb{B}_n^0(X)$.
Hence we obtain the following theorem.

\begin{theorem}
\label{thm:natural-Hn+1-module-strucrure-on-EnX}
There is a natural action of the Hecke algebra
$\mathcal{H}_{n+1}$ on $E_n^0(X)$
for any space $X$ such that
${\rm inc}: E_n^0(X)\to \mathbb{B}_n^0(X)$
is a map of $\mathcal{H}_{n+1}$-modules.
\end{theorem}

\section{Comparison of Hecke operators}
\label{section:comparison-Hecke-operators}


In this section we compare the Hecke operators
in $E_n$-theory with those in $E_{n+1}$-theory 
via $\mathbb{B}_n$-theory.
In particular,
we construct a ring homomorphism
$\omega: \mathcal{H}_{n+1}\to\mathcal{H}_n$
between Hecke algebras
and show that 
the $\mathcal{H}_{n+1}$-module structure
on $E_n^0(X)$ given 
by Theorem~\ref{thm:natural-Hn+1-module-strucrure-on-EnX}
is obtained 
from the $\mathcal{H}_n$-module
structure 
by the restriction along $\omega$
(Theorem~\ref{thm:comparison-Hn-Hn+1-on-EnX}).

First, we compare Hecke operators
in $\mathbb{B}_n$-theory with those
in $E_{n+1}$-theory.
By Proposition~\ref{prop:on-Phi-M-ch-and-equivariance},
we obtain the following proposition.

\begin{proposition}
We have 
$\widetilde{\rm T}_M^{\mathbb{B}}\circ {\rm ch}=
   {\rm ch}\circ \widetilde{\rm T}_M^{E_{n+1}}$
for any finite abelian $p$-group $M$
with {\rm $p$-rank$(M)\le n+1$}.
\end{proposition}

Next, we consider a relationship between
$\widetilde{\rm T}_M^{\mathbb{B}}$ and 
$\widetilde{\rm T}_N^{E_n}$.
In order to describe the relationship,
we will introduce a set
$I_n(M,N)$ and an integer
$a_n(M,N)$ for finite abelian $p$-groups $M$ and $N$,
where $p$-rank$(M)\le n+1$ and
$p$-rank$(N)\le n$.
We take a split surjection
$\mu: \Lambda^{n+1}\to \mathbb{Q}_p/\mathbb{Z}_p$,
and set $V=\ker\mu$.
We define $I_n(M,N)$ to be the set of all  
subgroups $M'$ of $\Lambda^{n+1}$
such that $M'\cong M$ and
$M'\cap V\cong N$:
\[ I_n(M,N)= \{M'\le \Lambda^{n+1}|\
    M'\cong M,\, M'\cap V\cong N\}.\]


\begin{lemma}
The cardinality $|m_n(N)|$ divides $|I_n(M,N)|$.
\end{lemma}

\proof
We fix an isomorphism $V\cong \Lambda^n$
and identify $V$ with $\Lambda^n$ via
this isomorphism.  
Let $f:I_n(M, N)\to m_n(N)$ be
a map given by $f(M')=M'\cap V$.
We take $N'\in m_n(N)$,
and suppose 
that $M'\in f^{-1}(N')$.
If $N''\in m_n(N)$,
then there is an automorphism $g\in {\rm Aut}(V)$
such that $g(N')=N''$. 
We can extend $g$ to an automorphism
$\widetilde{g}\in {\rm Aut}(\Lambda^{n+1})$
such that $\mu\circ \widetilde{g}=\mu$ 
and $\widetilde{g}|_V=g$.
We see that $\widetilde{g}(M')\in I_n(M,N)$ and
$f(\widetilde{g}(M'))=N''$.
Hence, if $f^{-1}(N')=\{M_1',\ldots,M_r'\}$,
then $f^{-1}(N'')=\{\widetilde{g}(M_1'),\ldots,\widetilde{g}(M_r')\}$.
This shows that $|m_n(N)|$ 
divides $|I_n(M,N)|$.
\qqq

We set
\[ a_n(M,N)=\displaystyle\frac{|I_n(M,N)|}{|m_n(N)|}.\]
We notice that $a_n(M,N)$ is independent of the choice
of a split surjection $\mu$.

By Theorem~\ref{thm:natural-Hn+1-module-strucrure-on-EnX},
the Hecke algebra $\mathcal{H}_{n+1}$
naturally acts on $E_n^0(X)$ for any space $X$.
We denote by 
$(\widetilde{\rm T}_M^{\mathbb{B}})^{\mathbb{G}_{n+1}}$
the operation associated to a finite abelian $p$-group $M$
with $p$-rank $\le n+1$.

\begin{proposition}\label{prop:relation-Hecke-TMB-to-TNEn}
We have 
\[ (\widetilde{\rm T}_M^{\mathbb{B}})^{\mathbb{G}_{n+1}}
    =\sum_{[N]}\, a_n(M,N)\ \widetilde{\rm T}_N^{E_n}, \]
where the sum ranges over all isomorphism classes 
of finite abelian $p$-groups with $p$-rank $\le n$.   
\end{proposition}

In order to prove Proposition~\ref{prop:relation-Hecke-TMB-to-TNEn},
we need the following lemma.
Let $i_{\mu}^{\mathbb{B}}:
\mathbb{B}_n^0\to D_{\mathbb{B}}(\mu)$ be
the $\mathbb{B}_n^0$-algebra structure map.

\begin{lemma}\label{lemma:i-mu-faithfully-flat}
The ring homomorphism $i_{\mu}^{\mathbb{B}}$ is faithfully flat.
\end{lemma}

\proof
By definition,
$D_{\mathbb{B}}(\mu)=\ \subrel{r}{\rm colim}
D(r,\mu[p^r])$.
The lemma follows from the fact that
$D(r,\mu[p^r])$
is a finitely generated free $\mathbb{B}_n^0$-module
for all $r$.
\qqq


\proof[Proof of Proposition~\ref{prop:relation-Hecke-TMB-to-TNEn}]
By Proposition~\ref{prop:relation-PhiB-PhiEn},
we have
\[ (i_{\mu}^{\mathbb{B}}\otimes 1)\circ
   \widetilde{\rm T}_M^{\mathbb{B}}\circ {\rm inc} =
   \sum_{M'\in m_{n+1}(M)} ({\rm inc}(\mu)\otimes {\rm inc})\circ
   \Phi_{M'\cap V}^{E_n}. \]
A decomposition
$m_{n+1}(M)=\coprod_{[N]}I_n(M,N)$
implies that the right hand side is 
\[ 
    ({\rm inc}(\mu)\otimes{\rm inc})\circ(i\otimes 1)\circ  
    \left(\sum_{[N]} a_n(M,N)\ \widetilde{\rm T}_N^{E_n}\right).
\]
Using 
$({\rm inc}(\mu)\otimes{\rm inc})\circ (i\otimes 1)=
(i_{\mu}^{\mathbb{B}}\otimes 1)\circ {\rm inc}$
and the fact that 
$i_{\mu}^{\mathbb{B}}$ is faithfully flat 
by Lemma~\ref{lemma:i-mu-faithfully-flat},
we obtain 
\[ \widetilde{\rm T}_M^{\mathbb{B}}\circ{\rm inc}=
   {\rm inc}\circ \sum_{[N]} a_n(M,N)\ \widetilde{\rm T}_N^{E_n},\]
which implies the desired formula.
\qqq

Conversely,
we shall write $\widetilde{\rm T}_N^{E_n}$
in terms of $(\widetilde{\rm T}_M^{\mathbb{B}})^{\mathbb{G}_{n+1}}$.
We define a partial order 
on the set of isomorphism classes of 
finite abelian $p$-groups of $p$-rank $\le n$
as follows.
We write $[A]\le [B]$
if there exists a monomorphism from $A$ to $B$.
We define $b_n(B,A)$ inductively as follows.
We set  $b_n(A,A)=1$.
When $[A]<[B]$, we set 
\[ b_n(B,A)= -\sum_{[A]\le [C]<[B]} a_n(B,C)\, b_n(C,A).\]
\if
For $[N]<[M]$,
we have
\[ \sum_{[N]\le [K]\le [M]}b_n(M,K)\, a_n(K,N)=0\]
by definition. 
Note that we also have
\[ \sum_{[N]\le [K]\le [M]} a_n(M,K)\, b_n(K,N)=0 \]
for $[N]<[M]$.
\fi

\begin{proposition}\label{prop:description-TEn-by-TBn}
For any finite abelian $p$-group $N$ 
with {\rm $p$-rank$(N)\le n$},
we have
\[ \widetilde{\rm T}_N^{E_n}=
  \sum_{[N']\le [N]} b_n(N,N')\,
   (\widetilde{\rm T}_{N'}^{\mathbb{B}})^{\mathbb{G}_{n+1}}.\]
\end{proposition}

\proof
We prove the proposition by induction 
on the partial ordering.
When $N=0$, the proposition obviously holds.
We suppose that the proposition holds for 
abelian $p$-groups $<[N]$.
By Proposition~\ref{prop:relation-Hecke-TMB-to-TNEn}, 
the right hand side is
\[ \widetilde{\rm T}_{N}^{E_n}
             + \sum_{[N']<[N]}a_n(N,N')\, \widetilde{\rm T}_{N'}^{E_n}
             + \sum_{[N']<[N]}b_n(N,N')\, 
               (\widetilde{\rm T}_{N'}^{\mathbb{B}_n})^{\mathbb{G}_{n+1}}.\]
By the hypothesis of induction,
we have
\[ \begin{array}{rcl}
    \displaystyle\sum_{[N']<[N]}a_n(N,N')\, 
    \widetilde{\rm T}_{N'}^{E_n}&
   = &\displaystyle\sum_{[N']<[N]}\sum_{[N'']\le [N']} 
      a_n(N,N')\, b_n(N',N'')\, 
     (\widetilde{\rm T}_{N''}^{\mathbb{B}_n})^{\mathbb{G}_{n+1}}\\[2mm]
  &= &\displaystyle\sum_{[N'']<[N]}
      \left(\sum_{[N'']\le [N']<[N]} a_n(N,N')\,
      b_n(N',N'')\right) 
     (\widetilde{\rm T}_{N''}^{\mathbb{B}_n})^{\mathbb{G}_{n+1}}.\\
   \end{array}  \]
Since $\sum_{[N'']\le [N']<[N]} a_n(N,N')\,
      b_n(N',N'')=-b_n(N,N'')$,
the right hand side is $\widetilde{\rm T}_{N}^{E_n}$.
\qqq

Now, 
we consider a relationship between
$(\widetilde{\rm T}_{\mathbb{B}}\left(p^r\right))^{\mathbb{G}_{n+1}}$ 
and 
$\widetilde{\rm T}_{E_n}\left(p^s\right)$,
where 
$(\widetilde{\rm T}_{\mathbb{B}}\left(p^r\right))^{\mathbb{G}_{n+1}}$
is the restriction
of the operation
$\widetilde{\rm T}_{\mathbb{B}}\left(p^r\right)$
to the $n$th Morava $E$-theory of spaces
through ${\rm inc}$. 

\begin{theorem}\label{thm:Hecke-B-to-En}
We have 
\[ (\widetilde{\rm T}_{\mathbb{B}}\left(p^r\right))^{\mathbb{G}_{n+1}}=
   \sum_{s=0}^r\ p^{(r-s)n}\, \widetilde{\rm T}_{E_n}\left(p^s\right).\]
\end{theorem}

\proof
We take a split surjection
$\mu: \Lambda^{n+1}\to \mathbb{Q}_p/\mathbb{Z}_p$,
and set $V=\ker\mu$.
Fix a subgroup $N$ of $V$ with order $p^s$.
Let 
$J_n(p^r,N)$ be the set of all subgroups $M$
of $\Lambda^{n+1}$ with order $p^r$
such that 
$M\cap V=N$.
By Proposition~\ref{prop:relation-PhiB-PhiEn},
it is sufficient to show that
$|J(p^r,N)|=p^{(r-s)n}$.

If $M\in J(p^r,N)$,
then 
the map $\mu$
induces a surjective map 
$\overline{\mu}: M\to \mathbb{Q}_p/\mathbb{Z}_p[p^{r-s}]$
with kernel $N$.
Let $a$ be a generator of 
$\mathbb{Q}_p/\mathbb{Z}_p[p^{r-s}]$.
If $x\in M$ such that $\overline{\mu}(x)=a$,
then $M$ is generated by $N$ and $x$.
Note that $p^{r-s}x\in\ker\overline{\mu}=N$.
Conversely, 
if we have $x\in\mu^{-1}(a)$ such that $p^{r-s}x\in N$,
then $\langle N,x\rangle\in J(p^r,N)$,
where $\langle N,x\rangle$ 
is the subgroup generated by $N$ and $x$.

We denote by $X(p^r,N)$ the set
of all $x\in\mu^{-1}(a)$ such that $p^{r-s}x\in N$.
By the above argument,
we see that $J(p^r,N)\cong X(p^r,N)/N$.
Let $N'$ be the subgroup of $V$
consisting of $z\in V$ such that $p^{r-s}z\in N$.
Since $X(p^r,N)$ is an $N'$-torsor,
we obtain 
$|J(p^r,N)|=|N'|/|N|=p^{(r-s)n}$.
\qqq

\begin{corollary}\label{cor:relation-TBb-TEn}
For $r>0$, we have
\[ \widetilde{\rm T}_{E_n}\left(p^r\right) =
   (\widetilde{\rm T}_{\mathbb{B}}\left(p^r\right))^{\mathbb{G}_{n+1}}
   -p^n (\widetilde{\rm T}_{\mathbb{B}}
   \left(p^{r-1}\right))^{\mathbb{G}_{n+1}}.\]
\end{corollary}

By Theorem~\ref{thm:natural-Hn+1-module-strucrure-on-EnX},
there is a natural action of the Hecke algebra
$\mathcal{H}_{n+1}$ on $E_n^0(X)$ for any space $X$.
Thus, 
we have the $\mathcal{H}_n$-module structure and
the $\mathcal{H}_{n+1}$-module structure on $E_n^0(X)$.
Finally,
we compare these two module structures.

For this purpose,
we define an additive map 
$\omega: \mathcal{H}_{n+1}\to \mathcal{H}_n$
by
\[ \omega(M)=
    \sum_{[N]} a_n(M,N) N.\]

\begin{proposition}
\label{lemma:multiplicative-property-omega}
The map $\omega$ is a surjective algebra homomorphism.
\end{proposition}

By Propositions~\ref{prop:relation-Hecke-TMB-to-TNEn}
and \ref{lemma:multiplicative-property-omega},
we obtain the following theorem.

\begin{theorem}
\label{thm:comparison-Hn-Hn+1-on-EnX}
The $\mathcal{H}_{n+1}$-module structure
on $E_n^0(X)$ is obtained from
the $\mathcal{H}_n$-module structure 
by the restriction along the algebra homomorphism $\omega$.
\end{theorem}


\proof[Proof of Proposition~\ref{lemma:multiplicative-property-omega}]
The surjectivity follows from
Proposition~\ref{prop:description-TEn-by-TBn}.
Thus,
it suffices to show that
\begin{align}\label{eq:proof-theta-algebra-hom}
  \sum_{[M_3]}c_{n+1}(M_1,M_2;M_3)\,
   a_n(M_3,N_3)
   =
   \sum_{[N_1],[N_2]}a_n(M_1,N_1)\, a_n(M_2,N_2)\,
   c_n(N_1,N_2;N_3)
\end{align}
for any finite abelian $p$-groups $M_1,M_2,N_3$
with $p$-rank$(M_i)\le n+1\ (i=1,2)$
and $p$-rank$(N_3)\le n$.

We fix an isomorphism 
$\Lambda^{n+1}/M_{1,i}\cong \Lambda^{n+1}$ 
for each $M_{1,i}\in m_{n+1}(M_1)$.
For $M_{2,j}\in m_{n+1}(M_2)$,
we let $L_{i,j}$ be the subgroup in $\Lambda^{n+1}$
which contains $M_{1,i}$ and 
satisfies $L_{i,j}/M_{1,i}\cong M_{2,j}$
under the isomorphism.
We have
\[ \{L_{i,j}\}=\coprod_{[M_3]} c_{n+1}(M_1,M_2;M_3)\,
   m_{n+1}(M_3) \]
as multisets.
This implies
\begin{align}\label{eq:comparison-Lij-1}
   \{L_{i,j}\cap V\}=
   \coprod_{[M_3],[N_3]} c_{n+1}(M_1,M_2;M_3)\,
   a_n(M_3,N_3)\, m_n(N_3).
\end{align}

We fix an isomorphism
$(\mathbb{Q}_p/\mathbb{Z}_p)/\mu(M_{1,i})\cong
\mathbb{Q}_p/\mathbb{Z}_p$.
We denote by $\mu'$ the composition
$\Lambda^{n+1}\cong \Lambda^{n+1}/M_{1,i}
\stackrel{\mu}{\to}
(\mathbb{Q}_p/\mathbb{Z}_p)/\mu(M_{1,i})\cong
\mathbb{Q}_p/\mathbb{Z}_p$,
and set $V'=\ker \mu'$.
Notice that there is an isomorphism
$V/(M_{1,i}\cap V)\cong V'$ and that
this isomorphism induces an isomorphism
\[ (L_{i,j}\cap V)/(M_{1,i}\cap V)\cong
    M_{2,j}\cap V'. \]
Hence we see that
\begin{align}\label{eq:comparison-Lij-2}
    \{L_{i,j}\cap V\}=\coprod_{[N_1],[N_2],[N_3]}
    a_n(M_1,N_1)\,a_n(M_2,N_2)\,c_n(N_1,N_2;N_3)\,
    m_n(N_3)
\end{align}
as multisets.
Comparing (\ref{eq:comparison-Lij-1})
and (\ref{eq:comparison-Lij-2}),
we obtain (\ref{eq:proof-theta-algebra-hom}).
\qqq


\begin{thebibliography}{99}

\bibitem{Ando1}
M. Ando.
Isogenies of formal group laws and power operations 
in the cohomology theories $E_n$.
Duke Math. J. 79 (1995), no. 2, 423--485. 

\bibitem{Ando2}
M. Ando. 
Power operations in elliptic cohomology and 
representations of loop groups.
Trans. Amer. Math. Soc. 352 (2000), no. 12, 5619--5666. 

\bibitem{AHS}
M. Ando, M. J. Hopkins and N. P. Strickland. 
The sigma orientation is an $H_\infty$ map. 
Amer. J. Math. 126 (2004), no. 2, 247--334. 

\bibitem{AMS}
M. Ando, J. Morava and H. Sadofsky. 
Completions of ${\mathbf Z}/(p)$-Tate cohomology of periodic spectra. 
Geom. Topol. 2 (1998), 145--174.



\bibitem{DHS}
E. S. Devinatz, M. J. Hopkins, and J. H. Smith.
Nilpotence and stable homotopy theory. I.
Ann. of Math. (2) 128 (1988), no. 2, 207--241. 

\bibitem{Ganter}
N. Ganter. 
Orbifold genera, product formulas and power operations. 
Adv. Math. 205 (2006), no. 1, 84--133.


\bibitem{Goerss-Hopkins}
P. G. Goerss and M. J. Hopkins. 
Moduli spaces of commutative ring spectra,  
Structured ring spectra. 
London Math. Soc. Lecture Note Ser., 315 (2004), 151--200.


\bibitem{Greenlees-Sadofsky}
J. P. C. Greenlees, and H. Sadofsky.
The Tate spectrum of $v_n$-periodic complex oriented theories.
Math. Z. 222 (1996), no. 3, 391--405. 

\bibitem{Greenlees-Strickland}
J. P. Greenlees and N. P. Strickland. 
Varieties and local cohomology for 
chromatic group cohomology rings.
Topology 38 (1999), no. 5, 1093--1139. 

\bibitem{Hopkins-Smith}
M. J. Hopkins, and J. H. Smith.
Nilpotence and stable homotopy theory. II.
Ann. of Math. (2) 148 (1998), no. 1, 1--49. 

\bibitem{HKR1}
M. J. Hopkins, N. J. Kuhn, and D. C. Ravenel.
Morava $K$-theories of classifying spaces and 
generalized characters for finite groups. 
Algebraic topology (San Feliu de Gu\'{i}xols, 1990), 
186--209,
Lecture Notes in Math., 1509, Springer, Berlin, 1992. 

\bibitem{HKR2}
M. J. Hopkins, N. J. Kuhn, and D. C. Ravenel.
Generalized group characters and complex oriented cohomology theories.
J. Amer. Math. Soc. 13 (2000), no. 3, 553--594. 

\bibitem{Hovey}
M. Hovey.
Bousfield localization functors and 
Hopkins' chromatic splitting conjecture. 
The \v{C}ech centennial (Boston, MA, 1993), 225--250,
Contemp. Math., 181, Amer. Math. Soc., Providence, RI, 1995. 

\bibitem{Hovey-Sadofsky}
M. Hovey, and H. Sadofsky.
Tate cohomology lowers chromatic Bousfield classes.
Proc. Amer. Math. Soc. 124 (1996), no. 11, 3579--3585. 



\bibitem{KM}
N. M.Katz and B. Mazur.
Arithmetic moduli of elliptic curves.
Annals of Mathematics Studies, 108. 
Princeton University Press, Princeton, NJ, 1985. 



\bibitem{Lurie}
J. Lurie,
Elliptic Cohomology III: Tempered Cohomology. 
available at 
https://www.math.ias.edu/\~{}lurie/.

\bibitem{MRW}
H. R. Miller, D. C. Ravenel, and W. S. Wilson.
Periodic phenomena in the Adams-Novikov spectral sequence.
Ann. of Math. (2) 106 (1977), no. 3, 469--516. 

\bibitem{Morava}
J. Morava.
Noetherian localisations of categories of cobordism comodules.
Ann. of Math. (2) 121 (1985), no. 1, 1--39. 

\bibitem{RZ}
M. Rapoport and Th. Zink.
Period spaces for p-divisible groups.
Annals of Mathematics Studies, 141. 
Princeton University Press, Princeton, NJ, 1996. 

\bibitem{Ravenel-green}
D. C. Ravenel.
Complex cobordism and stable homotopy groups of spheres.
Pure and Applied Mathematics, 121. 
Academic Press, Inc., Orlando, FL, 1986. 

\bibitem{Ravenel-red}
D. C. Ravenel.
Nilpotence and periodicity in stable homotopy theory.
Annals of Mathematics Studies, 128. 
Princeton University Press, Princeton, NJ, 1992. 

\bibitem{Rezk}
C. Rezk.
The units of a ring spectrum and a logarithmic cohomology operation.
J. Amer. Math. Soc. 19 (2006), no. 4, 969--1014. 

\bibitem{Rezk2}
C. Rezk.
The congruence criterion for power operations in Morava E-theory. 
Homology, Homotopy Appl. 11 (2009), no. 2, 327--379. 

\bibitem{Shimura}
G. Shimura.
Introduction to the arithmetic theory of automorphic functions.
Kan\^{o} Memorial Lectures, No. 1. 
Publications of the Mathematical Society of Japan, No. 11. 
Iwanami Shoten Publishers, Tokyo; 
Princeton University Press, Princeton, N.J., 1971. 

\bibitem{Stapleton1}
N. Stapleton.
Transchromatic generalized character maps. 
Algebr. Geom. Topol. 13 (2013), no. 1, 171--203. 

\bibitem{Stapleton2}
N. Stapleton.
Transchromatic twisted character maps. 
J. Homotopy Relat. Struct. 10 (2015), no. 1, 29--61. 

\bibitem{Strickland}
N. P. Strickland.
Finite subgroups of formal groups.
J. Pure Appl. Algebra 121 (1997), no. 2, 161--208. 


\bibitem{Torii0}
T. Torii.
The geometric fixed point spectrum of 
$(\mathbf{Z}/p)^k$ Borel cohomology for 
$E_n$ and its completion. 
Recent progress in homotopy theory (Baltimore, MD, 2000), 343--369,
Contemp. Math., 293, Amer. Math. Soc., Providence, RI, 2002. 


\bibitem{Torii1} 
T. Torii.  
On degeneration of formal group laws and
application to stable homotopy theory.
Amer. J. Math. 125 (2003), 1037--1077.


\bibitem{Torii3}
T. Torii.
Comparison of Morava $E$-theories. 
Math. Z. 266 (2010), no. 4, 933--951. 


\bibitem{Torii5}
T. Torii.
On $E_{\infty}$-structure of the generalized Chern character. 
Bull. Lond. Math. Soc. 42 (2010), no. 4, 680-690. 

\bibitem{Torii6}
T. Torii.
HKR characters, $p$-divisible groups
and the generalized Chern character.
Trans. Amer. Math. Soc. 362 (2010), no. 11, 6159--6181. 

\bibitem{Torii8}
T. Torii.
Comparison of power operations in Morava E-theories. 
Homology Homotopy Appl. 19 (2017), no. 1, 59--87. 
\end{thebibliography}
\end{document}